%%AMSTeX file
\documentstyle{amsppt}

\magnification=\magstep1
\overfullrule=0pt

\def\pageheight#1{\vsize#1\relax}   %%These two lines use the
\pageheight{47pc}                   %%measurements from AMSTeX version
\voffset-19pt                       %%2.1c, which are more convenient
                                    %%at \magstep 1 than version 2.1h.

%I added these new definitions:
\def\RTS{{{R_i(T,S)}}}
\def\oRT{{\overline{R_i(T,S)}}}
\def\oRTp{{\overline{R_i(T',S')}}}
\def\RDT{{{R_i(\delta', T,S)}}}

\def\YTp{{Y(T',S')}}
\def\YTpp{{Y(T'',S'')}}
\def\T{{(T,S)}}
\def\Tp{{(T',S')}}
\def\Tpp{{(T'',S'')}}
\def\hb{{\hfill\break}}
\def\oF{{\overline {F}}}

\def\oRjDT{\overline{R_{i,j}(\delta', T,S)}}
\def\oRDT{\overline{R_i(\delta', T,S)}}
\def\RDTp{{R_i(\delta', T',S')}}

\def\fap{{\frak{a}_P}}
\def\fapq{{\frak{a}_P^Q}}
\def\ofa{{\overline \frak{a}}}
\def\ofap{{\overline{\frak a_P}}}
\def\tF{{\widetilde F}}

\def\0{\emptyset}
\def\1{{-1}}
\def\g{{\gamma}}
\def\Ad{\mathop {\text{Ad}}\nolimits}
\def\Gqa{{G(\bQ ) \\ G(\bA )}}
\def\Hom{\mathop {\text {Hom}}\nolimits}
\def\Pga{{P(\bQ ) \\ G(\bA )}} 

\def\gv{{\sum_{\g \in V(\bQ)}}}
\def\hgv{{\sum_{\hg \in \hV(\bQ)}}}

\def\\{{\backslash}}
\def\bA{{\Bbb A}}
\def\bC{{\Bbb C}}
\def\bQ{{\Bbb Q}}
\def\bR{{\Bbb R}}
\def\bS{{\bold S}}
\def\bZ{{\Bbb Z}}
\def\cB{{\Cal B}}
\def\C{{\Cal C}}  
\def\cC{{\Cal C}}

\def\cS{{\Cal S}}

\def\dist{\mathop {\text{dist}}\nolimits}
\def\fG{{\frak g}}
\def\fS{{\frak s}}
\def\fa{{\frak a}}
\def\frac#1#2{{#1\over #2}}

\def\hD{{\widehat \Delta}}
\def\hf{{\hat f}}
\def\hg{{\hat \gamma}}
\def\htau{{\hat \tau}}
\def\hV{{\hat V}}
\def\oa{{\overline a}}
\def\oA{{\overline A}}
\def\oo{{\overline 0}}

\def\oH{{\overline H}}

\def\oR{{\overline R}}

\def\oX{{\overline X}}

\def\of{{\overline f}}

\def\supset{\mathrel{\supseteq}}
\let\mathsin=\sin
\def\sin{\mathrel{\subseteq}} 
\def\snot{\mathrel{\subsetneq}} 
\def\spn{\mathop {\text{span}}\nolimits}
\def\supn{\mathrel{\mathop{\text{sup}}\limits_{n\in \omega_N}}}
\def\tom{{\varpi}} 

\def\volume{\mathop {\text{vol}}\nolimits}
\def\vare{{\varepsilon}}
\def\vpi{{\tilde \pi}}

\def\fo{{\frak o}}  
\def\hfo{{\tilde{\frak o}}}  
\def\fO{{\frak O}}  
\def\hfO{{\tilde{\frak O}}}

\def\supp{\mathop {\text{supp}}\nolimits}
\def\ofa{\overline{\frak a}}

\def\oa{{\overline a}}
\def\oR{{\overline R}}

\def\SL{{\bold S \bold L}}
\def\obQ{\overline{\Bbb Q}}

\def\tpi{{\tilde \pi}}

\def\GL{{\bold G\bold L}}
\def\cvx{\mathop {\text{cvx}}\nolimits}
\def\sgn{\mathop {\text{sgn}}\nolimits}
\def\affspan{\mathop {\text{affspan}}\nolimits}
\def\interior{\mathop {\text{Int}}\nolimits}
\def\D{{\Delta}}
\def\cH{{\Cal H}}

\def\fb{{\frak b}}
\def\threenorm{{\|\mskip - 1.45mu|}}
\def\cR{{\Cal R}}
\def\b0{\text{\bf 0}}
\def\cB{{\Cal B}}

{\catcode`r=\catcode`@ \catcode`@=\catcode`c \def\noamslogo{\def\logo@{\relax}}}

\nologo

\topmatter

\title A Truncated Integral of the Poisson Summation Formula \endtitle

\author Jason Levy \endauthor

\address Department of Mathematics and Statistics, University of Ottawa,
585 King Edward, Ottawa, ON K1N 6N5 \endaddress

\email jlevy\@science.uottawa.ca \endemail

\thanks Partially supported by NSERC and NSF. \endthanks

\subjclass Primary 11F99, 11F72  \endsubjclass

\abstract
Let $G$ be a reductive algebraic group defined over $\bQ$, with
anisotropic centre.  Given a rational action of $G$ on a finite-dimensional
vector space $V$, we analyze the truncated integral of the theta series
corresponding to a Schwartz-Bruhat function on $V(\bA)$.  The Poisson
summation formula then yields an identity of distributions on $V(\bA)$. 
The truncation used is due to Arthur.
\endabstract
\endtopmatter

\document

\beginsection 0. Introduction

This paper is an extension of the previous paper [8].  In these two papers
we extend some of the results about integrating the Poisson summation
formula in Weil's famous paper [13], using methods
developed by Arthur to deal with infinities in the trace formula.
We will use many results from the paper [8], and that paper provides good
preparation for the complicated geometric constructions in section~3
and for the analysis in section 4.

Suppose that $V$ is a finite-dimensional vector space and that $f$ is
a Schwartz function on $V(\bA)$, the rational adelic points of $V.$
The Poisson summation formula tells us that
$$
\gv f(\g) = \hgv \hf(\hg),
$$
where $\hV$ is the vector space dual to $V$ and $\hf$ is the Fourier
transform of $f.$ Now suppose that a reductive algebraic group $G$
defined over 
$\bQ$ acts on $V$ via a rational representation 
$\pi: G \rightarrow \GL(V).$ The Poisson summation formula can now be
used to show that for every element $g$ of $G(\bA)$, 
$$
\gv f(\pi(g^\1) \g) = |\det \pi(g)| \hgv \hf (\tpi(g^\1) \hg),
\eqno (0.1)
$$
with $\tpi: G \rightarrow \GL(\hV)$ the representation of $G$
contragredient to $\pi$. Define the function $\phi_{f,\pi}$ on 
$G(\bA)$ by 
$$
\phi_{f,\pi}(g) = \gv f(\pi(g^\1)\g);
$$
then the equality (0.1) gives a relation between the two functions
$\phi_{f,\pi}$ and $\phi_{\hf,\tpi}$ on $G(\bA).$  Notice that for any
$f$ and $\pi$, the function $\phi_{f,\pi}$ is left $G(\bQ)$-invariant.

Weil [13] noticed that when the dimension of the representation is small 
compared to the rank of the group, the function $\phi_{f,\pi}$ is
integrable.  When this occurs, we may integrate the function
$\phi_{f,\pi}$ over $\Gqa$ and produce an equality of two
$G(\bA)$-invariant distributions through the integration of the
equality (0.1).   One then easily obtains (see [8])
an equality of sums of orbital integrals on $V$ and on $\hV.$
This equality is the basis of the Siegel-Weil formula.  For general
representations, however, the function $\phi_{f,\pi}$ will not be
integrable over $\Gqa$ and a more refined approach must be used.  We must
somehow regularize the functions $\phi_{f,\pi}$ and $\phi_{\hf,\tpi}$
before integrating.

The regularization that we use is a truncation invented by Arthur, and
described in [8].  The important properties of this truncation are reviewed
in the next section.  The problem is then to determine the
behaviour of the integral of the truncation of $\phi_{f,\pi}$ as a function
of the truncation parameter $T.$   In the previous paper we dealt with
reductive algebraic groups with rational rank at most two and with
anisotropic centre.  In this paper we deal with arbitrary reductive
algebraic groups with anisotropic centre. Given a general reductive
algebraic group $G$, the algebraic subgroup $G^1$ defined as the kernel of
all rational characters of $G$ is reductive with anisotropic centre, and
$G$ is the product of
$G^1$ and the maximal split torus $Z$ in the centre of $G$.  The action of
$Z$ can be introduced by making a Shintani zeta function, but we will not
discuss this further in this paper.  See [14] for more about Shintani zeta
functions.

We prove that the positive Weyl chamber is a finite union of sub-cones
depending only of $\pi$ satisfying the following property: within each
sub-cone, the integral of the truncation of $\phi_{f,\pi}$ with
respect to the point $T$ asymptotically approaches the value of a
finite sum of products of polynomials in $T$ and exponentials of
linear functionals of $T$. These functions do depend on the sub-cone. 

Given a sub-cone $\C$, write $J_\C (f,\pi)$ for the constant term of the
analytic function of the previous paragraph (see section 4 for our
definition of constant term). 
Then the basic form of the truncated Poisson
summation formula is the statement
$$
J_\C (f,\pi) = J_\C (\hf,\tpi).
$$
We can do better than this, as in the Selberg-Arthur trace formula, by
breaking up both sides as sums over certain equivalence classes in
the respective vector spaces; this is the way we present the material.

Notice that we obtain several truncated formulas this way.
Because our proof does not explicitly produce the numbers
$J_\C (f,\pi)$, it is not clear whether they are actually distinct for
different sub-cones $\C.$ 

The truncated Poisson summation formula developed here is potentially
useful both for producing a 
``regularized'' Siegel-Weil theorem, extending results from [13],
(see [7] for a different approach to this), and for new results about
Shintani zeta functions.

The results in this paper were obtained during postdoctoral fellowships at
the Institute for Advanced Study and Oklahoma State University.  The
author also wishes to thank IHES, MSRI, and the University of Toronto,
for their hospitality during the preparation of this paper. The author
thanks J. Bernstein, M. McConnell and D.\ Witte for helpful conversations.

\beginsection 1. Preliminaries.

Let $G$ be a reductive algebraic group over $\bQ$ with anisotropic centre.
Let $P_0$ be a minimal parabolic subgroup 
of $G$, and $M$ a Levi component and $N$ the unipotent radical of $P_0$,
all defined over $\bQ.$  These will be fixed throughout the paper.
Write $A$ for the maximal $\bQ$-split torus of $G$
contained in $M$ and $\fa$ for the corresponding real vector space
$\Hom(X(A)_\bQ,\bR).$  We will use the
phrase ``parabolic subgroup'' to denote a standard parabolic subgroup,
that is, a parabolic subgroup $P$ of $G$ defined over $\bQ$ and
containing $P_0.$ 
Given a parabolic subgroup $P$ we can define $N_P$, its unipotent radical,
$M_P$ the unique Levi component of $P$ containing $M$, and $A_P$,
the split component of the centre of $M_P.$    We
will use standard notation (see [8] or [1]) for roots, weights, and so
forth.  In particular, $\fa^+$ denotes the points in $\fa$ where all
positive roots are positive. 

Fix a maximal compact subgroup $K=\prod_v K_v$ of $G(\bA)$ that is
admissible relative to $M$, as in [2].
With this choice we can define as in [1] a continuous map
$H:G(\bA) \rightarrow \fa$  invariant under right
multiplication by $K.$  Given a parabolic subgroup $P\sin G$, write
$H_P$ for the projection of $H$ to $\fa_P.$  Then $H_P$ is left
$P(\bQ)$-invariant, and so can be seen as a map from $\Pga$ to
$\fa_P.$

Given $P\sin Q$ parabolic subgroups of $G$ and $X$ a point or subset of
$\fa$ [resp. $\fa^*$], write $X_P^Q$ for the projection of $X$ to
the subspace $\fa_P^Q$ of $\fa$ [resp. ${\fa^*}_P^Q$ of $\fa^*$], and
$\exp(X_P^Q)$ for the pre-image of $X_P^Q$ in $A_P^Q(\bR)^0$ under the
map $H_P$. To
simplify notation, we will remove all occurrences of $P_0$ as a subscript;
for example the Levi component of $P_0$ is $M=M_{P_0}.$  Notice that since
the centre of $G$ is anisotropic, $\fa_G=\b0$ and so we can also
eliminate all occurrences of $G$ as a superscript.  When it will cause no
confusion we may also simplify expressions with subscripts and superscripts
such as $P_1$ as follows: we will write $\fa_1^2$, for example, for
$\fa_{P_1}^{P_2}$ where $P_1 \sin P_2 $ are parabolic subgroups.

In the previous paper [8] we proved a result relating certain functions
on quotients of $G(\bA)$ and certain functions on $\fa.$  This result is
the key geometric principle that will allow us to produce the truncated
Poisson summation formula.

We first recall the functions involved.  We fix throughout this paper
a point $T_1$ in $-\fa^+$ and a compact subset $\omega$ of
$N(\bA)M(\bA)^1$ such that 
$$
G(\bA) = P(\bQ) \{pak \mid p\in \omega,\ a\in A(\bR)^0,\ k\in K,\
 \alpha(H(a) - T_1) > 0 \ \text{for all } \alpha\in \Delta^P\}
$$
for every parabolic subgroup $P\sin G.$  Then for a given $T\in \fa^+$, write
$F^P(\cdot,T)$ for the characteristic function of the relatively compact set
of all points $g\in P(\bQ)\\G(\bA)$ with a representative in the set
$$
\eqalign{
\{pak\mid p\in \omega,\ a\in A(\bR)^0,\ k\in K,\ &
\alpha(H(a)-T_1) > 0\  \text{for all } \alpha\in\Delta^P, \cr
&\tom(H(a) -T) \le 0\ \text{for all } \tom\in\hD^P\}.\cr}
$$
Given $T\in\fa^+$, we also define $\Gamma(\cdot,T)$ by
$$
\Gamma(X,T)=\sum_{R\colon P\sin R\sin Q} (-1)^{\dim(A_R / A_Q )} 
\tau_P^R (X) \htau_R^Q (T-X),
\quad X\in\fa,
$$
and for parabolic subgroups $P\sin Q$ we define $\Gamma_P^Q(\cdot,T)$ by
$\Gamma_P^Q(X,T)=\Gamma(X_P^Q,T_P^Q),$ $X\in\fa.$
Lemma 3.4 of [8] states that the function $X\mapsto \Gamma_P^Q(X,T),\,
X\in\fa_P^Q$, is the characteristic function of a convex set whose
closure is the convex hull of the points $T_R,\ P\sin R\sin Q.$

A relation between the functions $F$ and $\Gamma$ is given in Corollary
3.3 of [8] and runs as follows:  Let $T_2\in\fa^+$ be a fixed sufficiently
regular point.  Then for $T\in T_2+\fa^+,\ S\in \fa^+$, and $g\in\Gqa$,
$$
\eqalignno{
&F^G (g,T+S) - F^G (g,T) & (1.1)\cr
& = \ \sum_{P\sin Q \snot G} \sum_{\delta \in P(\bQ )
\\ G(\bQ )} F^P (\delta g,T_2)
\Gamma_P^Q  (H_P (\delta g) - T_2,T-T_2 )
\Gamma_Q (H_Q (\delta g) -T,S) \cr}
$$

Let $\pi$ be a rational representation of $G$ on a finite-dimensional
vector space $V.$  The following definitions were introduced in [8] and
require only that $G$ be a reductive algebraic group defined over $\bQ$.

Define a {\it semisimple vector} in $V$ to be one whose geometric orbit
$\pi(G(\obQ))\g \i V(\obQ)$ is
Zariski closed, and define a {\it nilpotent vector} in $V$ to be one such
that the origin is contained in the Zariski closure of its geometric
orbit.  If the representation $\pi$ is the Adjoint representation of $G$ on
its Lie algebra $\fG$, then these definitions give the
standard notions of semisimple and nilpotent elements of 
$\fG.$  Since the properties of being semisimple and nilpotent are both
invariant under the action of $G(\obQ)$, we will call an orbit, either
rational or geometric, {\it semisimple} or {\it nilpotent} if its elements
are.

Let $\gamma$ be an element of $V(\bQ ).$  A standard argument in invariant
theory shows that the Zariski closure of the geometric
orbit $\pi(G(\obQ))\g$ of $\g$ contains a unique closed geometric
orbit in $V(\obQ).$  The set of rational points in this closed
$G(\obQ)$-orbit is a union of semisimple $G(\bQ)$-orbits in
$V(\bQ)$---a non-empty union by the rational Hilbert-Mumford Theorem
([6]).    We define the {\it geometric semisimple component} of the
element $\gamma\in V(\bQ)$ to be this union.    If the representation $\pi$
is the Adjoint representation then the geometric semisimple component of an
arbitrary element  $\g\in V(\bQ)$ is the set of rational points in the
geometric orbit of the semisimple part $\g_s$ of $\g$ given by the Jordan
decomposition,  which may be larger than the $G(\bQ)$-orbit of $\g_s.$
Define a {\it geometric equivalence class} to be the set of
elements in $V(\bQ)$ whose geometric semisimple component equals a given
closed geometric orbit. 

\remark{Remark} For general representations, there is no known natural
assignment of a single semisimple vector to each vector, extending the
assignment to an element of its Lie algebra to its semisimple component in
the Jordan decomposition, and such an assignment is unlikely to exist.
(See [12] for an analysis of the case of
direct sums of the Adjoint representation.)  However, it may be possible to
identify a canonical rational (rather than geometric, as done above)
semisimple orbit to a given vector, by the following procedure, and that
would be sufficient for our purposes.
Define the Hilbert-Mumford closure of an orbit $\pi(G(\bQ))\g,\ \g\in
V(\bQ)$, to be the set of points in $V(\bQ)$ that can be expressed as
$$
\lim_{t\rightarrow 0} \pi(p(t)g) \g
$$
for some homomorphism $p:{\Bbb G}_m \rightarrow G$ defined over $\bQ$ and
some $g$ in $G(\bQ).$
Luna's property (A) (see [10]) asserts that the topological closure of every
orbit $\pi(G(\bR))v$ of a point $v\in V(\bR)$ contains a unique
$G(\bR)$ orbit of semisimple elements.  Together with the rational
Hilbert-Mumford theorem ([6]), property (A) implies that when the base
field is $\bR$ the Hilbert-Mumford closure of the rational orbit of a vector
contains a unique rational  semisimple orbit.  This latter property
trivially holds over algebraically closed fields, and is also known to hold
over $p$-adic fields [11], but appears to be unknown in general over
$\bQ.$  Some progress towards this question can be found in [9].
If it were  
true, then we would define the {\it semisimple component} of $\g$ to be
this unique semisimple orbit, and the {\it equivalence class} of $\g$ to be
the set of elements in $V(\bQ)$ whose semisimple component equals the
semisimple component of $\g.$  We note that if the
Hilbert-Mumford closure of every $G(\bQ)$-orbit contains a unique
semisimple orbit, then all the following results hold with equivalence classes
rather than geometric equivalence classes.
\endremark

The key property of the geometric equivalence classes is  the following
lemma, which generalizes Lemma 2.1 of [8] and is easily proven.

\proclaim{Lemma 1.1}  Let $G$ be a reductive algebraic group defined
over $\bQ$, let $\pi$ be a rational representation of $G$ on a
finite-dimensional vector space $V$, and let $A$ be a torus in $G$ that is
split over $\bQ.$  Suppose that $\Lambda_0$ and $\Lambda_+$ are two sets
of rational characters on $A$ such that there exists a point $a\in A(\bQ)$
satisfying 
$$
\eqalign{
|\lambda (a)|  =  1 & \quad\text{for every } \lambda \in \Lambda_0,\cr
|\lambda (a)|  >  1 & \quad\text{for every } \lambda \in \Lambda_+.\cr}
$$
Write
$$
V_0=\bigoplus_{\lambda\in\Lambda_0} V^\lambda, \quad 
V_+=\bigoplus_{\lambda\in\Lambda_+} V^\lambda,
$$
with $V^\lambda$ the weight space in $V$ corresponding to $\lambda.$ Then
for every geometric equivalence class $\fo$ in $V(\bQ)$ and subset $\cS$
of $V_0(\bQ)$, 
$$
(S + V_+ (\bQ )) \cap \fo \ = \ (\cS \cap \fo )
+ V_+ (\bQ ). 
$$
\endproclaim

Recall that the weight space in $V$ corresponding to 
$\lambda\in X^*(A)_\bQ$ is the vector subspace  $\{\,v\in V\mid
\pi(a)v=\lambda(a)v \hbox{\rm \ for all } a\in A\,\}$ of $V$.

Let $\fo$ be a geometric equivalence class in $V$, and define the function
$$
\phi_{\pi,\fo} (g,f) =  \sum_{\gamma \in \fo} \ f(\pi (g^{-1} ) \gamma ),
\quad g\in G(\bA),
$$
a left $G(\bQ)$-invariant function of $G(\bA).$  We will eliminate the
subscript $\pi$ when the representation is understood.  The Poisson
summation formula implies that
$$
\sum _{\fo \in \fO} \phi_{\pi,\fo} (g,f ) = \sum_{\hfo \in \hfO}
\phi_{\vpi,\hfo} (g,\hf), \eqno (1.2)
$$
where $\fO$ and $\hfO$ are the collections of geometric equivalence 
classes with respect to the representations $\pi$ and $\vpi$,
respectively---recall that we are assuming that the centre of $G$ is
anisotropic, so $|\det\pi(g)|=1$.

Given a Schwartz-Bruhat function $f$ on $V(\bA)$, a point $T\in\fa^+$, and
a geometric equivalence class $\fo\in\fO$, define
$$
J_\fo^T(f,\pi) = \int_{\Gqa} F^G(g,T) \phi_{\pi,\fo} (g,f ) dg.
$$
As in [8] we will try to find asymptotic formulas for $J_\fo^T$ for $T$ in
certain cones (to be defined later) in $\fa^+$ by examining differences
$J_\fo^{T+S}(f,\pi) - J_\fo^T(f,\pi)$ where $T$ and $S$ both lie in one
cone and $\|S\|\le 1.$  As on page 1392--3 of [8], we can write
$J_\fo(T+S)-J_\fo(T)$ as the finite sum over pairs $(P,Q)$ of
parabolic subgroups of $G$ with $P\sin Q\snot G$ of
$$
\eqalignno{
\int_{N_P(\bQ)\\N_P(\bA)} \int_{A_P(\bR)^0}& e^{-2\rho_P(H_P(a))} 
\phi_\fo (na,f^{M_P,K,T_2})  &(1.3)\cr
&\times \Gamma_P^Q(H_P(a),T-T_2) 
\Gamma_Q(H_P(a) -T,S) dadn,\cr}
$$
where
$$
f^{M_P,K,T_2}(v) = e^{-2\rho_P((T_2)_P^Q)} \int_{\omega_{M_P}} 
F^P(m,T_2) \int_K f(\pi(\exp((T_2)_P^Q)mk)^\1 v) dk dm,
$$
with $\omega_{M_P}$ a fundamental domain of $M_P(\bQ ) \\ M_P(\bA )^1$
and $T_2\in\fa^+$ a fixed sufficiently regular point.  Estimating each
of the integrals (1.3) requires some lengthy preliminary constructions
on the vector space $V$; these constructions of the subject of section 3.

\beginsection 2. A result of Brion-Vergne

In this section, we recall a result ([4], Theorem 4.2) on Fourier
transforms of convex polyhedra and give some consequences.  The proof
of our main theorem relies on the application of this result to certain
polytopes to be constructed in section 3.  Before we can state
Brion-Vergne's result, we must introduce some notation.

Suppose that $V$ is a finite-dimensional real vector space.  In
section 3, we will take $V=\fa_P$, so will write $V^*$ instead of
$\hV$ for the dual of $V$.  Suppose
that $\mu_1,\ldots,\mu_N$ are elements of $V^*$.  Given 
$x=(x_1,\ldots,x_N)\in\bR^N$ we
can define the convex polyhedron
$$
P(x) := \{\, v\in V\mid \mu_i(v) +x_i \ge 0, \ 1\le i\le N\,\}.
$$
Assume that for some $x\in \bR^N$, the polyhedron $P(x)$ is non-empty
and contains no line; this implies that $\mu_1,\ldots,\mu_N$ span
$V^*$.  Write $\cB$ for the set of subsets $\sigma$ of
$\{\,1,\ldots,N\,\}$ such that $\{\,\mu_i\mid i \not\in \sigma\,\}$ is
a basis of $V^*$.  Given $\sigma\in\cB$ we define the following three
objects:
$\{u_{i,\sigma}\}_{i\not\in\sigma}$ is the basis of $V$ dual to the
basis  $\{\,\mu_i\mid i \not\in \sigma\,\}$ of $V^*$,
$s_\sigma\colon \bR^N\to V$ is the linear map sending 
$x=(x_1,\ldots ,x_N)\in\bR^N$ to the 
unique point $v\in V$ with $\mu_i(v)+x_i=0$, $i\not\in \sigma$, and
$C(\sigma)$ is the set $\{\,x\in \bR^N\mid s_\sigma(x)\in P(x)\,\}$.

One can verify using an alternative description of $C(\sigma)$---the
one given in [4] 3.1, 4.1---that $C:=\cup_{\sigma\in\cB} C(\sigma)$
is the set of $x\in\bR^N$ such that $P(x)$ is non-empty.  To any set
$\Sigma\i \cB$, write $C_\Sigma$ for the intersection
$\cap_{\sigma\in\Sigma} C(\sigma)$, and to any $x\in C$ assign the set
$\Sigma_x=\{\,\sigma\in\cB\mid x\in C(\sigma)\,\}$.  The each set
$C_{\Sigma_x}$, $x\in C$, is a closed convex cone containing $x$, and
we obtain finitely many cones this way.  If we set 
$\cB(\g):=\{\,\sigma\in\cB\mid \g\i C(\sigma)\,\}$ for any such cone
$\g$, then $\g=C_{\cB(\g)}.$  If $\g$ is maximal among these cones, we
call it a {\it chamber} of C (in [4] it would be called the closure of a
chamber).

The following result is Theorem 4.2 of [4].

\proclaim{Theorem 2.1} Suppose that $V$ is a finite-dimensional real
vector space, that $\mu_1,\ldots,\mu_N$ are in $V^*$, and that at some
point in $\bR^N$, (2.1) defines a non-empty polyhedron containing no line.
Suppose that $\g$ is a chamber in $C$ and that $x$ is a point in
$\g$.\hb
(a) The extreme points of $P(x)$ are the $s_\sigma(x)$,
$\sigma\in\cB(\g)$, with possible repetition.\hb
(b) For generic
$$
\mu\in\{\,\sum_{i=1}^N c_i \mu_i \mid c_i \ge 0\,\},
$$
we have the identity
$$
\int_{P(x)} e^{-\mu(v)}dv =
\sum_{\sigma\in \cB(\g)} e^{-\mu(s_\sigma(x))}
{\volume \{\,\sum_{i\not\in\sigma} t_i u_{i,\sigma}\mid t_i\in[0,1]\,\}
\over \prod_{i\not\in\sigma}\mu(u_{i,\sigma})},
\eqno(2.2)
$$
where $dv$ and $\volume$ denote the same Lebesgue measure on $V$.
\endproclaim

\remark{Remark:} Notice that if (2.2) holds for a given $\mu$, then
the integral, as a function of $x\in\bR^N$, is a finite linear
combination of exponentials of linear functionals.
\endremark

We will need to tweak somewhat the statement of this theorem to suit
our needs.

\proclaim{Lemma 2.2} With the assumptions in Theorem 2.1, if for any
$x\in C$, the set $P(x)$ is bounded, then the set
$$\{\,\sum_{i=1}^N c_i \mu_i \mid c_i \ge 0\,\}$$
is all of $V^*$.
\endproclaim

\demo{Proof}  This follows from the theory of the polar.  It follows
immediately from [5] 6.1(a) and 9.1(b).\qed\enddemo

Recall that the bounded polyhedra are just the polytopes, the convex
hulls of finite sets of points.

Recall from [8] that we call a function $f$ on $\bR^N$ {\it t-finite}
if it is expressible (uniquely) as a finite sum
$\sum p_\lambda(x)e^{\lambda(x)}$, with $\lambda\in (\bR^N)^*$,
$p_\lambda\in\bC[x_1,\ldots,x_N]$.

\proclaim{Lemma 2.3} Keep the assumptions of Lemma 2.2.  Then for any
$\mu\in V^*$, the function
$$
x \mapsto F_\mu(x):=\int_{P(x)} e^{\mu(v)} dv,\quad x\in \g
$$
is t-finite.  In fact, each polynomial $p_\lambda$ in its
decomposition has degree at most $n$ ($n-1$ for $\lambda\ne 0$).
\endproclaim

\demo{Proof} Let $\mu\in V^*$, and pick $\mu_0\in V^*$ so that for
every $t\in(0,1)$, $\mu+t\mu_0$ is generic.  For any $x\in\g$
consider the function $t\in\bR \mapsto F_{\mu+t\mu_0}(x)$.  It is
a continuous function whose value at $t=0$ is $F_\mu(x)$.  On the
other hand, for $t\in(0,1)$, its value is given by (2.2).  Expanding
the exponentials in the right-hand side of (2.2) into a power series
in $t$ and letting $t\to 0$ gives the desired result. \qed\enddemo

\proclaim{Lemma 2.4}  Keep the assumptions of Lemma 2.2, and let
$y=(y_1,\ldots,y_N)$ be a point in $\bR^N$.  If for each extreme point
$v_x$ of $P(x)$, the set $\{\,v\in V\mid \text{if } \mu_i(v_x)+x_i=0,
\text{ then } \mu_i(v)+y_i=0,\ i=1,\ldots,N\,\}$ contains a single point
$v_y$ that is extreme in $P(y)$, then $y\in \g$.
\endproclaim

\demo{Proof} Fix $\sigma\in\cB(\g)$.  By Theorem 2.1, $s_\sigma(x)$ is
an extreme point of $P(x)$.  The extra hypothesis of the Lemma implies
that $s_\sigma(y)$ is an extreme point of $P(y)$, in particular
that $s_\sigma(y)\in P(y)$, so that $y\in C(\sigma).$  Therefore
$y\in C_{\cB(\g)}=\g$.\qed\enddemo

\remark{Remarks} (1) The map from the extreme points of $P(x)$ to
those of $P(y)$ given in the Lemma is surjective, by Theorem
2.1(a).\hb
(2) There is a surjection from the power set of $\{\,1,\ldots,N\,\}$
to the set of faces of $P(x)$, $x\in \g$, sending a set $S$ to the
set of points $v\in P(x)$ satisfying
$$
\mu_i(v)+x_i =0, \quad \hbox{\rm for\ all\ } i\in S, \eqno(2.3)
$$
or equivalently, to the convex hull of those extreme points of
$P(x)$ satisfying (2.3).  Under the hypotheses of Lemma 2.4, this
leads to a surjection from the set of faces of $P(x)$ to the set of
faces of $P(y)$.  If the hypotheses of the Lemma also hold
with $x$ and $y$ reversed, then this map between faces is a bijection.
\endremark

\beginsection 3. Geometry on $\fa.$

As mentioned in the introduction, we will find that the integral
$J_\fo^T(f,\pi)$ does not have simple asymptotic behaviour as $T$
approaches infinity in $\fa^+$, but rather that there exist convex open
cones in $\fa^+$ such that the function $J_\fo^T(f,\pi)$ is asymptotic to a
t-finite function in $T$ as $T$ tends to infinity within each
cone. Let us first describe these cones.

Let $\Psi_\pi \i \fa^*$ be the union of the set $\Delta$ of simple
roots of $(G,A)$ and the set of non-zero weights of $\pi$ 
with respect to $A.$  Define a function $d$ on $\fa$ as follows: for 
$X\in\fa,\ d(X)$ is the minimum, over all parabolic subgroups $P\sin Q
\snot G$ and subsets $\cS \sin \Psi_\pi$ such that $\spn(\cS_P)\cap
\fa_Q^*$ does not equal the trivial subspace $\b0,$ of the distance 
$$
\dist(\ker \cS,\cvx (X_R)_{P\sin R \sin Q})
$$
between the kernel $\ker \cS\i \fa$ of the
set $\cS$ and the convex hull 
$\cvx (X_R)_{P\sin R \sin Q} \i \fa_P \sin \fa$
of the projections of $X$ to $\fa_R$, $P\sin R \sin Q.$  Notice that for
$t\in \bR^+$ and $X\in \fa$, $d(tX)=td(X).$  This function $d$ is the
correct analogue to this situation of the functions denoted similarly in
[3] and [8].

Let us examine this definition through some examples.  Suppose that the
rational rank of $G$ is two.  If the set $\cS$ contains two
non-proportional elements (for example if $\cS = \Delta$), then for every
$P\sin Q \snot G$ the intersection $\spn(\cS_P)\cap \fa_Q^*$ is not $\b0$,
and 
$$
\dist(\ker \cS,\cvx (X_R)_{P\sin R \sin Q}) = \|X_Q\|.
$$
Suppose on the other hand that $\spn(\cS)$ is a line.  If this line does
not equal
$\fa^{1*}$ for a parabolic subgroup $P_1$, then the set of pairs $(P,Q)$ of
parabolic subgroups $P\sin Q \snot G$ with $\spn(\cS_P)\cap \fa_Q^* \ne
\b0$ is $\{(P,P) \mid P \ne G\}$.
If $\spn(\cS)$ does equal $\fa^{1*}$, then the set of pairs $(P,Q)$ is just
$\{(P,P) \mid P \ne P_1,\; G\}$, because the projection $\cS_1$ of $\cS$
to $\fa_1^*$ is zero.  In either case,
$$
\dist(\ker \cS,\cvx (X_R)_{P\sin R \sin Q}) = \dist(\ker \cS, X_P).
$$

%%%Pairs $(P,Q)$ with $P\ne Q$ may appear in groups of rank at least three.
%%%For example
%%%if $\Delta=\{\alpha_1, \alpha_2, \alpha_3\}$, $P=P_0$,
%%%$\Delta^Q=\{\alpha_3\}$, and 
%%%$\cS = \{\alpha_3 - \alpha_1, \alpha_3 -\alpha_2\}$, then
%%%$\dist(\ker \cS,\cvx (X_R)_{P\sin R \sin Q})$ measures the distance
%%%between
%%%a line and a line segment in $\fa \cong \bR^3$.  Pairs $(P,Q)$ with $P_0
%%%\snot P \snot Q$ first appear in rank four.  For example, 
More typical is the following example, with $G=\SL(5)$ and $\pi$ the
irreducible representation of $G$ whose highest weight is three times the
sum of the fundamental weights.  Here
$\Delta=\{\alpha_1, \alpha_2, \alpha_3, \alpha_4\}$.  Choose $P$, $Q$, and
$\cS$ by  $\Delta^P=\{\alpha_4\}$,
$\Delta^Q=\{\alpha_3,\alpha_4\}$, and $\cS = \{\alpha_4 + \alpha_3
-\alpha_1, -\alpha_4 + \alpha_3 -\alpha_2\}$. Then 
$\spn(\cS)\cap \fa_Q^* = \b0$, but $\spn(\cS_P)\cap \fa_Q^*$ is a line
containing the projection $(\alpha_2-\alpha_1)_Q= -\alpha_1 + \alpha_2 +2/3
\ 
\alpha_3 + 1/3 \ \alpha_4$ of $\alpha_2-\alpha_1$ to $\fa_Q^*$.  Here
$\dist(\ker \cS,\cvx (X_R)_{P\sin R \sin Q})$ measures the distance
between the plane $\ker \cS$  and the line segment $\cvx(T_P,T_Q)$ in
$\fa\cong \bR^4$, which may 
be less than the distance between the line $\ker(S)\cap\fa_P =
\ker\{\alpha_4, \alpha_3-\alpha_1, \alpha_3 -\alpha_2\}$ and
$\cvx(T_P,T_Q)$.

The set of zeros of $d$ is a finite union of closed convex cones, each
corresponding
to a single choice of $\cS$, $P$, and $Q$.  Each of these cones has
positive codimension; this is because $\ker \cS \cap \cvx(X_R)_{P\sin R
\sin Q}$ is contained in $\fa_P$ and hence in 
$\ker (\cS_P) \cap (X_Q + \fa_P^Q)$ and this set is non-empty only if $X_Q$
lies in $\ker(\spn (\cS_P) \cap \fa_Q^*)$, a non-trivial condition by the
non-triviality of $\spn(\cS_P)\cap\fa_Q^*$.  If we choose $P=Q$ and
$\cS$ to be a singleton $\{\lambda\}$, the subspace 
$\ker(\spn (\cS_P) \cap \fa_Q^*)$ is the hyperplane $\ker \lambda_P$,
so $d$ maps each of these hyperplanes to 0.  We will see that the
function $J_\fo^T(f,\pi)$ is asymptotic to a real-analytic function as
$T$ tends to infinity in each convex open cone in the complement in
$\fa^+$ of the set of zeros of $d$. 

Notice that the above facts about the zeros of $d$ trivially imply
that the complement to $d=0$ in $\fa^+$ can be expressed as a finite union
of convex open cones that may intersect.
Call the cones appearing in such a decomposition the $\pi$-dependent cones
in $\fa^+.$  The asymptotics of the function  $J_\fo^T(f,\pi)$ as $T$ varies
in the different $\pi$-dependent cones in $\fa^+$ will in general be
different; it seems likely however that their constant terms (when defined
correctly, see the comments before Theorem 4.5) will often be the
same.  If $\pi$ is the Adjoint representation, then $d$ is the
function given in [3], and there is only one $\pi$-dependent cone,
namely the whole of $\fa^+.$ 

\remark{Remark}  The complement of the set of zeros of $d$ in $\fa^+$
decomposes as a finite union of {\it disjoint} convex open cones in $\fa^+$.
We will not require this fact in the paper, so we will not prove it here. 
\endremark

Choose a non-empty convex open cone in the complement in $\fa^+$ of the set
of zeros of $d$ and call it $\C.$  The cone $\C$ will be fixed
throughout this section.  Given $\vare>0$, we write $\C_\vare$ for the set
of $X$ in $\C$ satisfying $d(X) > \vare\|X\|.$  Fix for the remainder of
this section $\vare > 0$ sufficiently small that $\C_\vare$ is a non-empty
convex open cone in $\fa^+.$  We also write $\C_\vare(1)$ for the set of
$X\in\C_\vare$ with $\|X\| \le 1.$  

The inner integral in (1.3) is over
the set of $a\in A_P(\bR)^0$ such that
$$
\Gamma_P^Q(H_P(a),T-T_2) \Gamma_Q(H_P(a) -T,S) =1,
$$
that is, $H_P(a)$ lies in a polytope in real Euclidean space $\fa_P.$  The
key to estimating (1.3) is breaking up this polytope into
pieces on which the integral is much easier to evaluate.

Fix for the remainder of this section two parabolic subgroups $P\sin Q
\snot G.$  Define $R_P^Q(T,S)$ as the support of the function
$X\mapsto \Gamma_P^Q(X,T) \Gamma_Q(X -T,S),$ $X\in\fa_P$, so that the integral
in (1.3) can be taken over $a$ with $H_P(a)\in R_P^Q(T-(T_2)_P^Q, S).$  We
will break 
up $R_P^Q(T,S)$ in a complicated way depending on the hyperplanes $\ker
\lambda,\ \lambda\in \Psi_{\pi,P}$ that intersect it.

For $T\in \fa^+$, define ${R'}_P^Q(T) \i T_Q + \fa_P^Q$ as the convex hull
$\cvx(T_R)_{P\sin R\sin Q}$ of the projections of $T$ to $\fa_R$; this is
the support in $T_Q+ \fa_P^Q$ of the function $X\mapsto\Gamma_P(X,T).$
Then for 
$T,S\in \fa^+$,
$$
R_P^Q(T,S) = {R'}_P^Q(T) + {R'}_Q(S) \supset {R'}_P^Q(T),
$$
so if $\|S\| \le 1$, every point in $R_P^Q(T,S)$ is within unit
distance of ${R'}_P^Q(T).$

We say that a hyperplane $H$ in an affine space is a boundary hyperplane of
a convex polytope with non-empty interior if $H$ is the affine span of a
codimension one face (a facet) of the polytope, and that a half-space is a
boundary half-space of a polytope it contains if its bounding hyperplane is
a boundary hyperplane.  Let the relative interior, the
relative boundary, and a relative boundary hyperplane of a region  be,
respectively, its
interior, its boundary, and a boundary hyperplane of the region considered
as an object {\it in its affine span.}
The dimension of a polytope means the dimension of its affine span.

\proclaim{Lemma 3.1}  For $T$ in $\fa^+$, the non-empty faces of
${R'}_P^Q(T)$ are exactly the convex hulls 
$$
\cvx(T_R)_{P_1\sin R \sin P_2} = {R'}_1^2(T),
$$
where $P_1$ and $P_2$ are parabolic subgroups satisfying $P\sin P_1 \sin
P_2 \sin Q.$  Furthermore, the dimension of ${R'}_1^2$ is $\dim \fa_1^2.$
\endproclaim

\demo{Proof} 
Notice that the affine span of ${R'}_P^Q$ is $T_Q+ \fa_P^Q.$  A face of
${R'}_P^Q$ in the intersection of ${R'}_P^Q$ with some of its relative
boundary hyperplanes.  Lemma 3.4 of [8] implies that the relative boundary
hyperplanes of ${R'}_P^Q$ are
$$
\{X \in T_Q + \fa^Q \mid \alpha(X) = 0\}, \eqno (3.1)
$$
for roots $\alpha\in\Delta_P^Q$, and
$$
\{X\in T_Q+ \fa_P^Q \mid \tom(X) = \tom(T)\}, \eqno(3.2)
$$
for weights $\tom\in\hD_P^Q.$
Consider a collection of these hyperplanes in $T_Q+ \fa_P^Q$ whose
intersection intersects ${R'}_P^Q$ non-trivially.  Choose parabolic
subgroups $P_1,P_2,\ P\sin P_1,P_2\sin Q,$ so that $\Delta_P^1$ is the set
of roots whose corresponding hyperplane (3.1) is in the collection and
$\hD_2^Q$ is the set of weights whose hyperplane (3.2) is in the collection.

Suppose that there were a root $\alpha\in \Delta_P^1$ whose dual weight
$\tom_\alpha$ in $(\fa_P^Q)^*$ lay in $\hD_2^Q.$  Since $\tom_\alpha$ can
be written in the form
$$
\tom_\alpha = c\alpha + \sum_{{\tom\in\hD_P^Q}\atop{\tom\ne \tom_\alpha}}
d_\tom \tom
$$
with $c$ positive and all $d$ nonnegative, we see that for every $X$ in the
intersection of the given collection of hyperplanes,
$$
\sum_{{\tom\in\hD_P^Q}\atop{\tom\ne \tom_\alpha}}
d_\tom \tom(X) = \tom_\alpha(X) = \tom_\alpha(T) = c\alpha(T) +
\sum_{{\tom\in\hD_P^Q}\atop{\tom\ne \tom_\alpha}} d_\tom \tom(T).
$$
Since $c\alpha(T)>0$ and all $d_\tom$ are non-negative, it must be that
$\tom(X)> \tom(T)$ for some $\tom\in\hD_P^Q.$  Therefore the intersection
of  ${R'}_P^Q$ with the collection of hyperplanes is trivial. This gives a
contradiction, so no such $\alpha$ can exist.  This proves that $P_1\sin
P_2.$

Therefore every non-empty face of  ${R'}_P^Q$ is the set of points $X$ in
$T_{P_2} + \fa_1^2$ such that 
$$
\eqalign{
\alpha(X) \ge 0 \quad&\text{for every } \alpha\in\Delta_P^Q, \cr
\tom(X) \le  \tom(T) \quad&\text{for every } \tom\in\hD_P^Q.\cr}
$$
This trivially equals the set of points $X$ in $T_{P_2} + \fa_1^2$ such
that
$$ 
\eqalign{
\alpha(X) \ge 0 \quad&\text{for every } \alpha\in\Delta_1^Q, \cr
\tom(X) \le  \tom(T) \quad&\text{for every } \tom\in\hD_P^2.\cr}
\eqno(3.3)
$$
Consider the smaller collection of inequalities 
$$
\eqalign{\alpha(X) \ge 0, \quad &\hbox{for every } \alpha\in\Delta_1^2 \cr
\tom(X)\le \tom(T), \quad &\hbox{for every } \tom\in\hD_1^2,\cr} \eqno
(3.4)
$$
for $X\in T_{P_2}+\fa_1^2$.
Since every weight $\tom\in\hD_P^2 \setminus \hD_1^2$ can be written in the
form
$$
\tom = \tom^1 + \sum_{\tom'\in\hD_1^2}d_{\tom'} \tom'
$$
with $\tom^1$ in $\hD_P^1$ and all $d_{\tom'}$ non-negative, if 
$X\in T_{P_2}+\fa_1^2$ satisfies the inequalities (3.4), it satisfies all
the inequalities in the second line of (3.3).  Since every root
$\alpha\in\Delta_1^Q \setminus\Delta_1^2$ can be written in the form 
$$
\alpha=\alpha_2 - \sum_{\tom\in\hD_1^2} d_\tom \tom 
$$
with $\alpha_2$ in $\D_2$ and all $d_\tom$ non-negative, if 
$X\in T_{P_2}+\fa_1^2$ satisfies the inequalities (3.4), then
$$
\eqalign{\alpha(X) &= \alpha_2(X) - \sum_{\tom\in\hD_1^2} d_\tom \tom(X)
\cr
&\ge \alpha_2(T) - \sum_{\tom\in\hD_1^2} d_\tom \tom(T) =\alpha(T) > 0,\cr}
$$
for any $\alpha\in\D_1^Q \setminus \D_1^2$, so that the inequalities (3.3)
and (3.4) are equivalent of $T_{P_2}+\fa_1^2$.
By Lemma 3.4 of [8], the set of $X\in T_{P_2} + \fa_1^2$ satisfying
(3.4) is exactly $\cvx (T_R)_{P_1\sin R\sin P_2}.$  This proves the first
statement of the Lemma.  The second statement is trivial, as $T$ lies in
$\fa^+.$  
\qed\enddemo

The following lemma points out the key properties of the region $\C$.

\proclaim{Lemma 3.2}  (i) The collection of subsets $\cS$ of $\Psi_{\pi,P}$
such
that $\ker \cS$ intersects ${R'}_P^Q(T)$ is independent of $T\in \C.$
\hfill\break
(ii) Given $T\in\C_\vare$, parabolic subgroups $P\sin Q \snot G$ and a
subset $\cS$ of $\Psi_{\pi,P}$, if the distance from $\ker \cS$ to
${R'}_P^Q(T)$ is less than $\vare \|T\|$, then
$\ker \cS$ must actually intersect ${R'}_P^Q(T)$.

\endproclaim

\demo{Proof} The following fact will be needed: 
Given $\cS \sin \Psi_{\pi,P}$, let $\cR \sin \Psi_\pi$ be such that $\cS 
=\cR_P;$ since the points of $\ker (\cR_P)$ closest to ${R'}_P^Q(T)\i
\fa_P$ actually lie in $\ker(\cR_P) \cap \fa_P \sin \ker(\cR)$, we
have the inequality
$$
\dist(\ker \cR,{R'}_P^Q(T)) \le \dist(\ker \cS,{R'}_P^Q(T)).
$$

(i). \quad Suppose that $\ker \cS$ intersects ${R'}_P^Q(T)$ but not
${R'}_P^Q(T')$, with $\cS\sin \Psi_{\pi,P},$ and $T,T'\in\C.$  Consider a
minimal face of ${R'}_P^Q(T)$ that intersects $\ker \cS$, so that the
intersection is a single point in the relative interior of the face. By
Lemma 3.1 this face is of the form
${R'}_1^2(T)$ for two parabolic subgroups $P_1,\,P_2,$ with $P\sin P_1\sin
P_2\sin Q$, and hence $\ker \cS$ intersects the affine space $T_{P_2} +
\fa_1^2=\affspan{R'}_1^2(T)$ in a point.  If we write $\cS_1$ for the set of
projections of weights in $\cS$ to $\fa_1^*$, this implies that 
$$
\spn(\cS_1, \fa_2^*) = \fa_1^*. \eqno(3.5)
$$

Now, by our assumptions, $\ker \cS$ intersects ${R'}_1^2(T)$ but not
${R'}_1^2(T').$  Therefore for some point $T''\in\C$ on the line segment
joining $T$ and $T'$, $\ker \cS$ intersects only the relative boundary of
${R'}_1^2(T'')$, not its relative interior. The intersection $(\ker\cS)
\cap {R'}_1^2(T'')$ is again a single point, and it must lie on a relative
boundary hyperplane of ${R'}_1^2(T'').$  First suppose that this relative
boundary hyperplane is of the form $(\ker \alpha_1) \cap (T_2+\fa_1^2)$,
for some $\alpha_1\in\Delta_1^2.$  By (3.5) we know that $\alpha_1$ is in
the
span of $\cS_1$ and $\fa_2^*$, so that 
$$
\spn(\cS_1 \cup \{\alpha_1\}) \cap \fa_2^* \ne 0.
$$
But this contradicts the facts that $\ker(\cS_1\cup\{\alpha_1\})$
intersects ${R'}_1^2(T'')$ and that $T''\in\C.$  Next suppose that the
relative boundary hyperplane of ${R'}_1^2(T'')$ is of the form 
$$\{X\in T_2+\fa_1^2 \mid \tom(X) = \tom(T)\}$$
for some $\tom\in\hD_1.$  Again, by
(3.5) we know that $\tom$ is in the span of $\cS_1$ and $\fa_2^*$, so if we
choose $P_3\snot P_2$ to be the parabolic subgroup with $\hD_3=\hD_2
\cup\{\tom\}$, we have
$$
\spn(\cS_1) \cap \fa_3^* \ne 0.
$$
But this contradicts the facts that $\ker\cS_1$ intersects ${R'}_1^3(T'')$
and that $T''\in\C.$  We have proven part (i).

(ii). \quad Pick $\cR$ as above, and consider a minimal face ${R'}_1^2(T)$
of ${R'}_P^Q(T)$ such that 
$$
\dist(\ker \cR,{R'}_1^2(T)) = \dist(\ker \cR,{R'}_P^Q(T)) \le \vare \|T\|.
$$
Then, since $T$ lies in $\C_\vare$, we know that $\spn(\cR_1) \cap \fa_2^*
= 0$, so that $\ker \cR$ intersects $T_1 + \fa_1^2$.  Clearly,
given a convex set $A$ and an affine subspace $S$ that intersects the
affine span of $A$ but not $A$ itself, the points on $A$ closest to $S$
lie on the relative boundary of $A$.  But then the minimality
of ${R'}_1^2(T)$ implies that $\ker \cR$ intersects ${R'}_1^2(T) \sin
{R'}_P^Q(T)$, which proves (ii).

%%%To prove the above statement, notice first that it can be easily
%%%demonstrated if $A$ and $S$ are zero- or one-dimensional, as the problem
%%%can then be situated in $\bR^2$.  In general, given a relative interior
%%%point $X$ of $A$, write $Y$ for the point on the intersection of $S$ and the
%%%affine span of $A$ that is closest to $X$, and $Z$ for the point of $S$
%%%closest to $X$.  Let $A'$ be the intersection of $A$ with the line joining
%%%$X$ and $Y$, and $S'$ be the minimum affine subspace containing $Y$ and $Z$
%%%(either a line or a single point), so that both $A'$ and $S'$ are at most
%%%one-dimensional, and hence there is a point on the relative boundary of
%%%$A'$ closer to $S'$ than is $X$.  This point is also on the relative
%%%boundary of $A$ and is closer to $S$ than is
%%%$X$, proving the statement.
\qed\enddemo

Given a functional $\lambda\in\fa_P^*$, write $H_\lambda$ for the
hyperplane $\ker \lambda$ of $\fa_P.$ Given in addition a positive
real number $b$, define a ``thickened hyperplane'' $H_\lambda(b)$ to
be the set of points $X\in\fa_P$ such that $|\lambda(X)| \le b.$ We
will say that $T\in\C$ is sufficiently large if $\|T\|$ is
sufficiently large.  We will also say $(T,S)$ is well-situated if
$S,T\in\C_\vare$, $T$ is sufficiently large, and $\|S\|\le 1$.  Recall
that $\vare>0$ is fixed throughout this section.

\remark{Remark} Since the largest angle between two vectors in the convex cone
$\C_\vare\i\fa^+$ is
less than $\pi/2$, it is valid to conclude that if $(T,S)$ and $(T',S')$
are both well-situated, then so is every point $(T'',S'')$ on the line
segment in $\fa\times\fa$ joining them.
\endremark

\proclaim{Lemma 3.3}  There exists $B\in\fa^*$, positive on $\C$, such that
for every parabolic subgroups $P\sin Q$, every subset $\cS$ of
$\Psi_{\pi,P}$, all well-situated, the intersection
$$
\bigcap_{\lambda\in\cS} H_\lambda(B(T)) \cap R_P^Q(T,S)
$$
is non-empty if and only if the intersection
$$
\bigcap_{\lambda\in\cS} H_\lambda \cap {R'}_P^Q(T) = (\ker \cS) \cap
{R'}_P^Q(T)
$$
is non-empty.
\endproclaim

\demo{Proof} We first observe that since $\Psi_{\pi,P}$ is finite 
there exists a positive constant $\kappa>1$ such that for any $b>0$ and
$\cS\sin \Psi_{\pi,P}$, if $X\in\fa_P$ lies in $\cap_{\lambda\in\cS}
H_\lambda(b)$, then $\dist(X,\ker \cS) < \kappa b.$  This follows from
repeated applications of the following: given a subspace $\bS$ of $\fa_P$
not contained in a hyperplane $H\i \fa_P$ through $0$, a point within
distance $B$ of both $\bS$ and $H$ is within $B/\mathsin(\theta/2)$ of
their intersection, where $\theta$ is the minimum angle between the
projections of $\bS$ and $H$ to $\fa_P / (\bS \cap H)$.
(Notice that the projection of $\bS$ is one-dimensional and not contained
in the projection of $H$, so the angle $\theta$ is  well-defined and
positive.  The distance from the given point to $\bS \cap H$ is the length
of the diagonal through $\bS\cap H$ of a parallelogram in $\fa_P/(\bS\cap
H)$ with one side in the projection of $\bS$, an adjacent side in the
projection of $H$, and all side lengths at most $B$.  This length is at
most $B/\mathsin(\theta/2)$, the length of the long diagonal of a rhombus
with inner angles $\theta \le \pi - \theta$ whose parallel sides are
distance $B$ apart.)

Fix a parabolic subgroup $Q\supseteq P.$  Since $R_P^Q(T,S)$ contains
${R'}_P^Q(T),\ \ker\cS \cap {R'}_P^Q(T) \ne \0 \ \hbox{implies }$
$$
\bigl( \bigcap_{\lambda\in
\cS}H_\lambda(b)\bigr) \cap R_P^Q(T,S) \ne \0,
$$
for any $b\ge0.$  On the other hand, if an intersection
$\cap_{\lambda\in\cS} H_\lambda(b)$ intersects $R_P^Q(T,S)$, then since
every point in $R_P^Q(T,S)$ is within unit distance of ${R'}_P^Q(T)$,
$$
\dist({R'}_P^Q(T), \bigcap_{\lambda\in\cS} H_\lambda ) \le
\dist(R_P^Q(T,S), \bigcap_{\lambda \in\cS} H_\lambda) +1 \le \kappa b +1.
$$
For $b\le {\vare \over 2\kappa} \|T\|$ and $T$ sufficiently large, this
implies that the distance between ${R'}_P^Q(T) = \cvx(T_R)_{P\sin R \sin
Q}$ and $\cap_{\lambda\in\cS} H_\lambda$ is less than $\vare \|T\|.$  Since
$T\in\C_\vare$, Lemma 3.2 lets us conclude that $\cap_{\lambda\in\cS}
H_\lambda$
intersects ${R'}_P^Q(T).$

We therefore need only find $B\in\fa^*$ such that $0< B(T)\le {\vare\over
2\kappa} \|T\|$ for all $T\in\C.$  The choice
$B= {\vare\over2\kappa\|\alpha\|}\alpha,$ $\alpha$ any simple root in
$\Delta$,
satisfies this condition.
\qed\enddemo

Let $B$ be as in Lemma 3.3, and let us return to the fixed parabolic
subgroups $P\sin Q \snot G$.  For the remainder of this section we
will work in $\fa_P$ so 
that by the kernel of a linear functional we will mean the kernel in
$\fa_P$, and by the interior of a region in  $\fa_P$, its interior in the
topology of $\fa_P$.  Write $\Pi=\Pi_P\in\fa_P^*$ for the
set of non-zero weights of $\pi$ with respect to the torus $A_P.$
Given any subset $\Pi^+$ of $\Pi$, write $R_P^Q(\Pi^+,T,S)$ for the closure
of the set of $X\in R_P^Q(T,S)\i \fa_P$ such that for $\lambda\in \Pi,\
\lambda(X)>0$ exactly when $\lambda \in \Pi^+.$  Clearly each non-empty set
$R_P^Q(\Pi^+,T,S)$ is a closed convex polytope with non-empty interior in
$\fa_P$, and for $(T,S)$ fixed, any 
two have disjoint interiors, and their (finite) union over $\Pi^+\sin \Pi$
is $R_P^Q(T,S).$ 

Let $\Pi^+\i \Pi$ be a set such that $R_P^Q(\Pi^+,T,S)$ is non-empty.
Given any subset $\Lambda\sin\Psi_{\pi,P}=\Pi\cup\Delta_P,$ define
$\Lambda^+ = \Lambda\cap\Pi^+,\ \Lambda^-=\Lambda\setminus \Pi^+,$ and
define the sign $\sgn (\lambda)$ of a weight $\lambda\in\Pi$ to be 1 if
$\lambda\in \Pi^+$ and $-1$ if $\lambda\in\Pi^-.$  For any functional
$\lambda \in \Delta_P \cup \Delta_Q \cup \hat{\Delta}_Q \cup
\hat{\Delta}^Q_P$, define $\sgn \lambda = 1.$  We will now decompose
$R_P^Q(\Pi^+,T,S)$ as a union of closed convex polytopes in $\fa_P$ with
disjoint non-empty interiors; in the following section we will evaluate the
contribution to (1.3) of the integral over each of these latter polytopes.
The constructions of these polytopes requires a lengthy recursion.

Let $\delta=1/|\Pi|$.  Then for any $\Pi^+\sin\Pi$, any positive
number $b$, any subset $\cS$ of $\Pi$, any linear combination
$$\mu=\sum_{\lambda\in\cS} d_\lambda \lambda,$$
and any $X\in\fa_P^+$ such that 
$$
\eqalign{(\sgn \lambda)\lambda(X) & >  0, \quad \hbox{for all }
\lambda\in\Pi,\cr
\mu(X) & > b, \cr}
$$
some $\lambda\in\cS$ with  $d_\lambda(\sgn \lambda) > 0$ must satisfy 
$$
\lambda(X) > \delta b/d_\lambda. \eqno (3.6)
$$

The $0^{\text{\it th}}$ step of the recursion proceeds as follows: for
each subset $\Lambda_0$ of $\Pi$, 
let $R_P^Q(\Lambda_0;\Pi^+,T,S)$ be the closure of the set of $X\in
R_P^Q(\Pi^+,T,S)$ such that for $\lambda\in\Pi$, $(\sgn\lambda)\lambda(X)
\ge B(T)$ exactly for $\lambda\in \Lambda_0.$  Notice that each non-empty
set $R_P^Q(\Lambda_0;\Pi^+,T,S)$ is a closed convex polytope in $\fa_P$,
that the non-empty sets corresponding to different $\Lambda_0$ have
disjoint nonempty interiors, and that their union over all
$\Lambda_0 \i \Pi$ equals $R_P^Q(\Pi^+,T,S).$

Now suppose that we have constructed, at the $k^{\text{\it th}}$ step,
a non-empty region 
$$R_P^Q(\Lambda_0,\ldots, \Lambda_k;\Pi^+,T,S),$$
with $\Lambda_0,\ldots,\Lambda_k$ disjoint
subsets of $\Pi.$  Our ``inductive hypothesis" is that the region
$R_P^Q(\Lambda_0,\ldots, \Lambda_k;\Pi^+,T,S)$ has non-empty interior, and
that  $X \in \fa_P$ belongs to $R_P^Q(\Lambda_0,\ldots,
\Lambda_k;\Pi^+,T,S)$ exactly when $X$ satisfies the following
inequalities:
$$
\eqalignno{
(\sgn\lambda)\lambda(X) \ge 0 \quad &\hbox{for all }
\lambda\in(\Pi\setminus \bigcup_{i=0}^k \Lambda_i) \cup \Delta_P, \cr
(\sgn\lambda)\lambda(X) \ge\delta_i B(T)\quad &\hbox{for all }
\lambda\in\Lambda_i,\ i=0,\ldots,k, & (3.7)\cr
\openup\jot
(\sgn\lambda)\lambda(X)\le\delta_i B(T)\quad &\hbox{for all }
\lambda\in\Lambda_{i+1},\ \hbox{if } i=0,\ldots,k-1,\cr
&\hbox{and all } \lambda\in\Pi\setminus\bigcup_{j=0}^k \Lambda_j,\ \hbox{
if
} i=k,\cr
\noalign{\noindent the inequalities of immediate concern to us, and}\cr
\tom(X) \le \tom(T), \quad &\hbox{for all } \tom\in\hD_P^Q,\cr
\alpha(X) \ge\alpha(T), \quad &\hbox{for all } \alpha\in\Delta_Q&(3.8)\cr
\tom(X) \le\tom(T+S),\quad &\hbox{for all } \tom\in\hD_Q.\cr}
$$
The constants $\delta_i$ are all positive, with $\delta_0=1.$ We will now
express $R_P^Q(\Lambda_0,\ldots, \Lambda_k;\Pi^+,T,S)$ as a union of
regions with disjoint non-empty interiors.

If the intersection
$$
\ker(\Pi \setminus \bigcup_{i=0}^k \Lambda_i) 
   \cap R_P^Q(\Lambda_0, \ldots, \Lambda_k; \Pi^+,T,S)
$$
is non-empty, we end the recursion and do not break up the set
$R_P^Q(\Lambda_0, \ldots, \Lambda_k;\Pi^+,T,S)$ any further.  Otherwise,
there exists an assignment $a$ of constants $c_a(\lambda)$ (independent of
$T,S$) for each functional $\lambda\in\cup_{i=0}^k\Lambda_i \cup \Delta_P$
such that the linear combination
$$
\mu_a= \sum_{\lambda\in\cup\Lambda_i \cup \Delta_P} c_a(\lambda) \lambda
\in\spn(\bigcup_{i=0}^k \Lambda_i\cup\Delta_P)
$$
lies in $\spn(\Pi\setminus\cup_{i=0}^k \Lambda_i)$, and the weighted sum of
certain of
the inequalities (3.7) is an inequality of the form $\mu_a(X) \ge c_a B(T)$
with the number $c_a$ strictly positive---this follows from a very special
case of the Krein-Milman theorem together with Lemmas 3.2 and 3.3.  Pick
one such assignment $a$ and write $\mu_a$ as a linear combination of
elements
of $\spn(\Pi\setminus\cup_{i=0}^k \Lambda_i)$:
$$
\mu_a = \sum_{\lambda\in\Pi\setminus\cup_{i=0}^k\Lambda_i} 
d_\lambda \lambda.
$$
Let $D$ be the maximum of $|d_\lambda|,\ \lambda\in
\Pi\setminus\cup_{i=0}^k\Lambda_i.$ By (3.6) we can conclude that for all
$X$ in $R_P^Q(\Lambda_0, \ldots, \Lambda_k;\Pi^+,T,S)$, at least one
$\lambda\in \Pi\setminus\cup_{i=0}^k\Lambda_i$ satisfies
$$
(\sgn \lambda)\lambda(X) \ge {\delta c_a\over D} B(T).
$$
Let $\delta_{k+1}$ equal $\delta c_a/D$, and for each subset $\Lambda_{k+1}$
of $\Pi\setminus\cup_{i=0}^k\Lambda_i$ define 
$$
R_P^Q(\Lambda_0, \ldots,\Lambda_{k+1};\Pi^+,T,S)
$$
to be the closure of the set of points $X$
satisfying the strict inequality of each inequality in (3.7) and (3.8) and
also
$$
\eqalign{
(\sgn \lambda)\lambda(X) &> \delta_{k+1} B(T) \quad  \hbox{for all }
\lambda\in\Lambda_{k+1} \cr
(\sgn \lambda)\lambda(X) &< \delta_{k+1} B(T) \quad  \hbox{for all }
\lambda\in\Pi\setminus\cup_{i=0}^{k+1} \Lambda_i.\cr}
$$
The non-empty sets $R(\Lambda_0, \ldots,\Lambda_{k+1};\Pi^+,T,S),$ are the
regions constructed at the $k+1^{\text{\it st}}$ step.  Notice that
these regions are determined by inequalities of the form (3.7) and
(3.8), that their union over all non-empty $\Lambda_{k+1}$ is 
$R_P^Q(\Lambda_0, \ldots, \Lambda_k;\Pi^+,T,S),$ and that their interiors
are pairwise disjoint.

Since $\Pi$ is finite and each $\Lambda_i, i\ge1$, is non-empty in a given
non-empty region 
$R_P^Q(\Lambda_0, \ldots, \Lambda_k;\Pi^+,T,S),$
we can write $R_P^Q(T,S)$ as a {\it finite} union
$$
\bigcup_{(\Lambda_0,\ldots,\Lambda_k;\Pi^+) \in I(T,S)}
R_P^Q(\Lambda_0, \ldots, \Lambda_k;\Pi^+,T,S) \eqno (3.9)
$$
of convex polytopes with disjoint non-empty interiors, indexed by a
(finite) set $I(T,S)$ of ordered tuples of varying size, such that the
region corresponding to each  tuple
$(\Lambda_0,\ldots,\Lambda_k;\Pi^+) \in I(T,S)$ is not broken up by the
above algorithm.  We will use this decomposition of $R_P^Q(T,S)$ in the
next section, by estimating an integral of the form of (1.3) over each of
these regions.  Notice that by Lemmas 3.2 and 3.3, the index set $I(T,S)$
and the implicit constants $\delta_i$ and other implicit choices made for
each element of $I(T,S)$, can be 
chosen independently of $S\in \C_\vare(1)$ and sufficiently large
$T\in\C_\vare.$  We can therefore denote this set simply as $I$; it depends
on $P$, $Q$, and, of course, $\pi.$  The set $I$ is in no way canonical,
but is sufficient for our purposes.

Consider a region $R_i(T,S)=R_P^Q(\Lambda_0, \ldots, \Lambda_k;\Pi^+,T,S)$
corresponding to 
$$i=(\Lambda_0,\ldots,\Lambda_k;\Pi^+) \in I(T,S),$$
and write 
$$
\Pi_0=\Pi\setminus\bigcup_{i=0}^k\Lambda_i, \quad 
\Pi_+ = \bigcup_{i=0}^k\Lambda_i^+ \sin \Pi^+, \quad
\Pi_- = \bigcup_{i=0}^k\Lambda_i^-\sin \Pi^-.
$$
Notice that by the
construction of $R_i(T,S)$, the set of weights of $\pi$ that vanish on
$\ker \Pi_0 \cap R_i(T,S)$ equals $\Pi_0$, so in particular $(\spn\Pi_0)
\cap\Pi= \Pi_0.$  The region $R_i(T,S)$ is a convex polytope bounded by the
inequalities (3.7) and (3.8), so that each boundary hyperplanes of
$R_i(T,S)$ is given by one of the following equations:
$$
\eqalignno{
\lambda(X) = 0, \quad &
\lambda\in(\Pi\setminus \bigcup_{i=0}^k \Lambda_i) \cup \Delta_P, \cr
(\sgn\lambda)\lambda(X) =\delta_i B(T),\quad & \lambda\in\Lambda_i,\
i=0,\ldots,k, & (3.7)'\cr 
\openup\jot
(\sgn\lambda)\lambda(X)=\delta_i B(T),\quad &
\lambda\in\Lambda_{i+1},\ \hbox{if } i=0,\ldots,k-1,\cr
&\hbox{and all } \lambda\in\Pi\setminus\bigcup_{j=0}^k \Lambda_j,\ \hbox{
if
} i=k,\cr
\noalign{\noindent and}\cr
\tom(X) = \tom(T), \quad & \tom\in\hD_P^Q,\cr
\alpha(X) =\alpha(T), \quad & \alpha\in\Delta_Q& (3.8)'\cr
\tom(X) =\tom(T+S),\quad & \tom\in\hD_Q.\cr}
$$

Let $Y$ be an extreme point of $R_i(T,S)$, let $\cH_1(T,S)$ and
$\cH_2(T,S)$ be the set of boundary hyperplanes of $R_i(T,S)$ of the form
$(3.7)'$ and $(3.8)'$, respectively, that contain $Y$.   Each of the
hyperplanes in $\cH_k(T,S),\ k=1,2$, depends explicitly on $T$ and $S$
through the corresponding equality in $(3.7)'$ or $(3.8)'$, so that  given
any $T', S' \in \fa$, we can naturally define the sets $\cH_k(T',S')$.

\proclaim{Lemma 3.4} Suppose that we are given $(T,S)$ well-situated,
$i\in I(T,S)$, and $Y$ an extreme point of $R_i(T,S)$.  Let $\cH_1$,
$\cH_2$ send elements of $\fa\times\fa$ to sets of hyperplanes in $\fa_P$
as described above.  If $(T',S')$ is also well-situated, then the
intersection
$$
\bigcap_{H\in\cH_1(T',S') \cup \cH_2(T',S')} H \eqno (3.10)
$$
has exactly one element.  Call this element $Y(T',S')$.  Then
$Y(T',S')$ is an extreme point of $R_i(T',S')$.
\endproclaim

%%%\remark{Remarks} (1) In other words, the extreme points of $R_i(T,S)$ are
%%%linear in $T$ and~$S.$ (Since we can interchange $T,S$ with $T',S'$ in the
%%%Lemma, it is clear that the extreme points of $R_i(T,S)$ and $R_i(T',S')$
%%%are in bijection.)
%%%
%%%(2) Lemma 3.4 and the following comment together imply that for
%%%$(T,S), (T',S')$ well-situated to every face $F$ of $R_i(T,S)$ corresponds
%%%a face $F'$ of $R_i(T',S')$ whose affine span is a translate of that of
%%%$F$, and intersects $\ker \lambda, \  \lambda \in \Pi_0$ arbitrary, if $F$
%%%does. 
%%%
%%%(3) Suppose
%%%we are given a collection of affine subspaces in a vector space, whose
%%%intersection is non-empty.  If each affine subspace in the collection is
%%%independently translated, and the intersection of all the translates is
%%%non-empty, then intersection of all the translates is a translate of the
%%%intersection of the original affine subspaces.

\demo{Proof} 
Let $\bS_{(k)}$ be $\cap_{H\in\cH_k(T,S)} H$, for $k=1,2$.  
Let $\cS\i \Psi_{\pi,P}$ be the set of functionals $\lambda$ such
that the point $Y$ satisfies an equality in $(3.7)'$ involving $\lambda$,
and let $R\sin Q$ be the parabolic subgroup containing  $P$ so that
$\hD_R^Q$ is the set of $\varpi\in\hD_P^Q$ whose hyperplane
$\varpi(X)=\varpi(T)$ appears in $\cH_2(T,S)$.  Then Lemma 3.3 and the
definition of $d$ imply that  $\spn (\cS) \cap \fa_R^* = \b0. $
Since we also know that $\bS_{(1)}\cap \bS_{(2)}=\{Y\}$, we can
conclude that 
$$
\fa_P^*=\spn\cS \oplus\fa_R^*, \eqno(3.11)
$$
and that the intersection of the hyperplanes $H_R$,
$H\in\cH_2(T,S)$, of $\fa_R$ is the point $Y_R$.  Since $Y$ lies in
$R_i(T,S)$, $Y_Q-T_Q$ must lie in $R_Q'(S)$, and projecting the
preceding sentence to $\fa_Q$, we find that the point 
$Y_Q-T_Q$ is extreme in $R_Q'(S)$.  Lemma 3.1 then says that for some
parabolic subgroup $P_1\supset Q$, $\{Y_Q-T_Q\}={R'}_1^1=\{S_1\}$, with
$S_1$ the projection of $S$ to $\fa_1$, so $Y_Q=T_Q+S_1$.  The definition
of $R$ tells us that $Y_R^Q=T_R^Q$, so that $Y_R=T_R+S_1$.  Putting all
this together, we see that $\bS_{(2)}=T_R+S_1+\fa_P^R$.

Going now to $T'$ and $S'$, notice that the hyperplanes in $\cH_1(T',S')$
are of all the form 
$$
\lambda(X)=kB(T'),
$$
where the constants $k$ and functionals $\lambda\in\Psi_{\pi,P}$ are
determined by $Y$.  The intersection $\cap_{H\in \cH_1(T',S')} H$ equals
${B(T')\over B(T)} \bS_{(1)}={B(T')\over B(T)}Y + \ker \cS$.  Also, the
previous paragraph implies the equality
$$
\bigcap_{H\in\cH_2(T',S')}H = T_R' + S_1' + \ker \fa_R^*. \eqno (3.12)
$$
We conclude from (3.11) that the set (3.10) contains a single point
$Y(T',S')$.  This proves the first statement of the Lemma.

Since $Y(T',S')$ lies on a collection of boundary hyperplanes of
$R_i(T',S')$ that intersect in a point, it suffices to prove simply that
$Y(T',S')$ lies in $R_i(T',S')$, that is, that $Y(T',S')$ satisfies all the
inequalities in (3.7) and (3.8).

We will first consider the inequalities in (3.8).  The inequalities on
the last two lines of (3.8) follow immediately from (3.12).  The
inequalities on the first line of (3.8) hold for $\varpi\in\hD_R^Q$ and are
actually equalities in the case, again by (3.12).  Lastly, let
$\varpi\in\hD_P^Q\setminus\hD_R^Q$.  Since $\fa_P^*=\spn\cS \oplus
\fa_R^*$, we know that 
$$\spn\cS\cap\fa_{R'}^*\ne\bold{0},\eqno(3.13)$$
where $R'\sin Q$
is the parabolic subgroup satisfying $\D_{R'}^Q=\D_R^Q\cup\{\varpi\}$, and
that $\varpi(Y)\le\varpi(T)$.  If $\varpi(Y(T',S'))>\varpi(T')$, then for
some $(T'',S'')$ on the line segment joining $(T,S)$ and $(T',S')$ (and by
an earlier remark, necessarily well-situated), the point $Y(T'',S'')$ would
satisfy $\varpi(Y(T'',S''))=\varpi(T'')$.  But then (3.13) and Lemma~3.3
would yield a contradiction.

We next show that $\YTp$ lies in $R_P^R(T',S')$.  Given the above
paragraph, this is equivalent to showing that $\alpha(\YTp)\ge0$ for all
$\alpha\in\D_P^R$.  
%%%Note: The $\D_R,\hD_R$ conditions are clear.  For $\hD_P^R$ conditions,
%%%let $\varpi\in\hD_P\setminus\hD_R$, so that $\varpi^R\in\hD_P^R$,
%%%$\varpi=\varpi^R+\varpi_R$.  Then $\varpi_R(\YTp)=\varpi_R(T),$$
%%%$\varpi(\YTp)\le\varpi(T)$, so get the desired ineq.
Assume otherwise, and let $\Tpp$ be the point on the
line segment joining $\T$ and $\Tp$ such that 
\item{$\bullet$} $\alpha(\YTpp)\ge0$ for all $\alpha\in\D_P^R,$
\item{$\bullet$} $\beta(\YTpp)=0$ for some $\beta\in\D_P^R$ with
$\beta(\YTp)<0$.

\noindent Then the point $\YTpp$ lies in
$$
\bigcap_{\lambda\in\cS\cup\{\beta\}} H_\lambda(B(T''))\cap R_P^R\Tpp,
$$
so by Lemma 3.3, the set $\ker(\cS\cup\{\beta\})\cap {R'}_P^R(T'')$ is
non-empty, and since $T''$ lies in $\cC_\vare$, we see that
$\spn(\cS\cup\{\beta\})\cap\fa_R^*=\bold{0}$.  By (3.11) this implies
that $\beta\in\spn\cS$.  But then $\beta(\YTp)$ and $\beta(\YTpp)=0$ can be
explicitly given as $cB(T')$ and $cB(T'')$, respectively, for some constant
$c$.  The constant $c$ must be zero, and we obtain a contradiction to the
definition of $\beta$.  Therefore $\YTp$ does in fact lie in $R_P^R(T',S')$.

Now, let $\mu$ be any linear functional in $\Psi_{\pi,P}$.  We must show
that the inequalities in (3.7) that mention $\mu$ hold for the point
$\YTp$ (and $T'$ instead of $T$).  If
$\mu\in\spn\cS$, then $\mu(Y)$ and $\mu(\YTp)$ are explicitly given  as
$cB(T)$ and $cB(T')$, respectively, for some constant $c.$  Since $Y$
satisfies all the inequalities in (3.7) (with $T$),  $\YTp$ must
satisfy those inequalities from (3.7) that mention $\mu$ (with
$T'$). Next, suppose that $\mu\not\in\spn\cS$.  Because of (3.11) we must
have that $\spn(\cS\cup\{\mu\})\cap\fa_R^*\ne\0$, so since $T\in\cC_\vare$,
$\dist(\ker(\cS\cup\{\mu\}),{R'}_P^R)>\vare\|T\|.$  Lemma 3.3 then says
that given a point $X\in R_P^R\T$, $|\lambda(X)|>B(T)$ for some
$\lambda\in\cS\cup\{\mu\}$. Now, $Y\in R_P^R\T$ and $|\lambda(Y)|\le B(T)$
for all $\lambda\in\cS$, so $(\sgn\mu)\mu(Y)>B(T)$.  The inequalities (3.7) 
must be consistent with this, so the only inequality in (3.7) that mentions
$\mu$ must be of the form $(\sgn\mu)\mu \ge \delta_k B(T)$, for some $k$. 
On the other hand, $\YTp\in R_P^R\Tp$, so that we can similarly obtain the
inequality
$(\sgn\mu)\mu(\YTp)>B(T')$.  Since $\delta_k\le1$, we conclude that $\YTp$
satisfies every inequality in (3.7) that mentioned $\mu$.  This finishes
our proof that $\YTp$ lies in, and hence is an extreme point of, $R_i\Tp$. 
\qed\enddemo

In the course of estimating an integral over $R_i(T,S)$ in the next
section, we will express the integrand as a sum of terms corresponding to
certain subsets of $\Pi_0$, and to estimate the integral of the piece of
the integrand corresponding to one subset $\Pi_1$ of $\Pi_0$, we will need
to further manipulate the set $R_i(T,S).$  The necessary constructions
form the remainder of this section.

Let $\Pi_1 \subseteq \Pi_0$ be a subset of $\Pi_0$ satisfying
$$
\Pi_1=\{\lambda \in \Pi_0|\lambda ((\ker \Pi_1) \cap R_i(T,S))=0\},
$$
so that in particular $(\spn \Pi_1) \cap \Pi_0=\Pi_1$.  (Note that both
$\{0\}\cap\Pi_0$ and $\Pi_0$ satisfy this property.)  Let $\cB \subseteq
\Pi_1$ be a basis and let $d$ be the dimension of $\spn \Pi_1$.  We want to
examine the dependence of the convex polytope 
$$
(X+ \ker \Pi_1)\cap R_i(T,S) \eqno (3.14)
$$
on $T, S$ and $X \in R_i (T,S)$.

Write $\overline{\fa_P}$ for the quotient space $\fa_P/\ker \Pi_1$, and
make the natural identification of the dual space of $\overline{\fa_P}$
with $\spn\Pi_1$. The projection map $\fa_P \to \overline{\fa_P}$ sending
$X \in \fa_P$ to its projection, to be denoted $\overline{X}$, in
$\overline{\fa_P}$ sends polytopes to polytopes, so the projection of
$R_i(T,S)$ is a polytope $\overline{R_i(T,S)}$.  

The dependence of the set (3.14) on $X$ is clearly through $\overline{X}$. 
To simplify matters, we will consider only points $\overline{X} \in
\overline{R_i(T,S)}$ close to zero, in a sense to be defined presently.

The point $\overline{0} =\ker \Pi_1$ lies in $\overline{R_i(T,S)}$ since
$\ker \Pi_1 \supset \ker \Pi_0$ intersects $\RTS$, and it is an extreme
point because of the inequalities $(\sgn \lambda) \lambda \geq 0, \ 
\lambda\in\Pi_1$, that hold on $\RTS$.  The polytope $\oRT$ has finitely
many facets $\oF$ through $\oo$; the boundary half-space corresponding to
one such facet $\oF$ is given by an inequality of the form
$$
\sum_{\lambda \in \cB} c_{\lambda}^\oF \lambda \geq 0, \eqno (3.15)
$$
for some numbers $c_{\lambda}^\oF,$ $\lambda\in\cB$.

\proclaim{Lemma 3.5} For any well-situated $\Tp$, the inequality (3.15)
defines a boundary half-space of $\oRTp$.
\endproclaim

\demo{Proof} That the half-space defined by (3.15) is a boundary half-space
of $\oRT$ is equivalent to the following statement:
The point $\sum_{\lambda\in\cB}c_{\lambda}^\oF \lambda $ belongs to and is
extreme among those $\mu\in\spn\Pi_1$ such that $\mu(X)\ge0$ for all
$X\in\RTS$.
We must prove that if we replace $\T$ with $\Tp$, this statement is still
true.

By the theory of the polar ([5], Theorem 6.4) and the extremality of 
$\sum_{\lambda\in\cB}c_{\lambda}^\oF \lambda $, the inequality (3.15)
can be written as a linear combination of the inequalities (3.7) and
(3.8).  Furthermore, if we take $X\in\RTS$ a 
pre-image of $\oo\in\oF$, then each inequality that appears in the above
linear combination with non-zero coefficient is actually an equality at
$X$.  Let $\Lambda$ [resp. $W$] be the set of functionals appearing in
equalities from $(3.7)'$ [resp. $(3.8)'$] that hold at $X$.  Then we have an
equality
$$
\sum_{\lambda\in\Lambda} d_\lambda\lambda +\sum_{\varpi\in W} d_\varpi\varpi
=
\sum_{\lambda \in \cB} c_{\lambda}^\oF \lambda,\quad
$$
for some constants $d_\lambda,\,d_\varpi$.  The set $\spn (W) + \fa_Q^*$ is
of the form $\fa_R^*$ for some parabolic subgroup $R\sin Q$.  The point $X$
lies in $R_P^R\T$ and in $\cap_{\lambda\in\Lambda\cup\cB} H_\lambda(B(T))$,
so by Lemma 3.3 and the well-situatedness of $\T$,
$\spn(\Lambda\cup\cB)\cap\fa_R^*=\0$, and so all the constants $d_\varpi$
must be 0.  Therefore the inequality (3.15) is a linear combination just
of inequalities from (3.7), and so takes the form
$$
\sum_{\lambda \in \cB} c_{\lambda}^\oF \lambda \geq cB(T),
$$
for some constant $c$.  The constant $c$ must clearly be 0.

Going now to $R_i\Tp$, we see that taking the same linear combination of the
inequalities (3.7) (with $T'$ instead of $T$), we find that the inequality
(3.15) holds also on $R_i\Tp$.  If the point
$\sum_{\lambda\in\cB}c_{\lambda}^\oF \lambda$ were a convex combination of
distinct functionals $\mu\in\spn\Pi_1$ that are non-negative on all of
$R_i\Tp$, then we could reason as above to find that these functionals
are also non-negative on all of $\RTS$, and so obtain a contradiction to
the above-stated extremality of $\sum c_\lambda^\oF \lambda$.  This
completes our proof.
\qed\enddemo

The functional 
$$
\lambda_{\cB}  = \sum_{\lambda\in \cB} (\sgn \lambda) \lambda
$$
is non-negative  on all $\RTS$ (and vanishes on the non-empty set
$(\ker \Pi_1) \cap \RTS$) for any well-situated $(T,S)$.  The proof of
Lemma 3.5 implies that the value of $\lambda_\cB$ at each extreme point
$Y(T,S)=\cap_{H\in \cH_1(T,S) \cup \cH_2(T,S)} H$ is determined from the
equalities in $(3.7)'$ that define hyperplanes $\cH_1(T,S)$, so that
$\lambda_\cB(Y(T,S)) = c_Y B(T)$ for a constant $c_Y$  independent of $T$
and $S$.  Therefore, by Lemma 3.4, for well-situated $(T,S)$, the minimal
non-zero value of 
$\lambda_\cB$ on the extreme points of $\RTS$ is given by $2\delta' B(T)$
for some non-zero constant $\delta'$ independent of $T$ and $S$.
Let $\oH$ be the hyperplane 
$\{X\in\ofa_P \mid \lambda_{\cB}(X)=\delta'B(T)\}$ of $\ofa_P$.  Since the
non-zero extreme points of $\oRT$ are projections of the extreme points
of $\RTS$ where $\lambda_{\cB}$ is non-zero, and $\lambda_{\cB}$ is left
invariant by projection to $\overline{\frak{a}}_P$, $\oH$  separates
$\oo$ from the other extreme points of $\oRT$.

Let $R_i(\delta', T,S)$ be the set of $X$ in $\RTS$ satisfying
$$
\lambda_{\cB}(X)=  \displaystyle \sum_{\lambda \in \cB} (\sgn \lambda)
\lambda (x) \leq \delta' B(T). \eqno (3.16)
$$
The boundary hyperplanes of $R_i(\delta', T,S)$ are exactly the boundary
hyperplanes of $\RTS$ that intersect $(\ker \Pi_1) \cap \RTS$, and the
hyperplane $\lambda_\cB=\delta'B(T)$ from (3.16).  This is because every
face of $\RTS$ that does not intersect $\ker \Pi_1$ is the convex hull of a
collection of extreme points of $\RTS$ not in $\ker \Pi_1$, and so has no
points satisfying (3.16).  This argument also
shows that the projection  $\overline{R_i(\delta', T,S)}$ of
$R_i(\delta', T,S)$ to $\ofap$ is a pyramid with apex $\oo$ whose boundary
half-spaces are exactly those given by (3.15) and (3.16).

We want to examine the dependence of (3.14) on $\oX, T,S$, with $(T,S)$
well-situated, and $\oX$ in the interior of $\overline{R_i(\delta', T,S)}$. 
A basic example is the intersection of translates of the line $y=z=0$
with the octahedron $R\i \bR^3$ given by
$$
0 \leq y+z, y-z, x+y-z, x+y+z \leq B(T).
$$
(the polytope $ 0 \leq y+z, y-z \leq \delta B(T), \ \delta B(T) \leq 
x+y-z, x+y+z \leq B(T)$ is similar); this models the points of the
intersection $(X+\ker \Pi_1)\cap \RTS$ in  $T_P+\fapq$ when $\dim \fapq =3$.
The intersection
$$
((x_0, y_0, z_0) + L) \cap R, \quad (x_0, y_0, z_0)\in \interior (R),
$$
is a line segment and each of its two endpoints lies on a boundary
hyperplane of $R$.  The two hyperplanes on which the end points lie are
$x+y-z=0, x+y+z=B(T)$ if $z_0 \geq 0$ and are $x+y+z=0, x+y-z=B(T)$ if $z_0
\leq 0$.   Therefore, the extreme points of 
$$
((x_0, y_0, z_0) + L) \cap R
$$
are linear in $T$ and in $(x_0, y_0, z_0)$ on a specified side of $z_0
= 0$.  We will see that the general situation is similar.

Notice that $\lambda_{\cB} (X+ \ker \Pi_1)= \lambda_{\cB}(X)$, so that for
any point $X$ in $R_i(\delta', T,S)$, the polytope (3.14) equals $(X+ \ker
\Pi_1) \cap R_i(\delta', T,S)$.

Given $X$ in the interior $\interior R_i(\delta', T,S)$ of $R_i(\delta',
T,S)$ and an extreme point $Y$ of the polytope (3.14), there is a face
$F$ of $\RDT$ such that 
$$
(X+ \ker \Pi_1) \cap \affspan F=\{\,Y\,\}. \eqno (3.17)
$$
Since $X\in\interior \RDT$, $\lambda_\cB(Y)<\delta'B(T)$, so $F$ is
not contained in the
hyperplane $\lambda_\cB=\delta'B(T)$. Therefore $F$ intersects $\ker\Pi_1$,
so that $\affspan F\cap \RTS$ is a face of $\RTS$ of dimension $\dim F$.
Since (3.17) contains a unique point, elementary linear algebra
implies that
$$
|(Y+\ker \Pi_1)\cap ( \affspan F)|\le1 \text{ for all } Y\in \fap.
\eqno (3.18) 
$$

Suppose that $(T', S')$ is also well-situated, that $F$ is a face of
$R_i(\delta', T,S)$ satisfying (3.18), and that $F'$ is the corresponding face
of $R_i(\delta', T',S')$ (more precisely, the intersection with
$R_i(\delta', T',S')$ of the face of $R_i(T',S')$ corresponding, by
Lemma 3.4 and the remark after Lemma 2.4, to $R_i(T,S)\cap\affspan F$;
$F'$ clearly also satisfies (3.18)).
Suppose also that $X\in \interior R_i(\delta', T,S), X' \in \interior
R_i(\delta', T',S')$ satisfy 
\item{$\bullet$} $(X+\ker \Pi_1)\cap \affspan F$ contains a (unique)
point, which
is extreme in $(X+\ker \Pi_1) \cap R_i(\delta', T,S)$,
\item{$\bullet$} $(X'+\ker \Pi_1)\cap \affspan F'$ contains a (unique)
point, which
is not extreme in $(X'+\ker \Pi_1) \cap R_i(\delta', T',S')$.

\noindent Then for some point $(X'', T'', S'') \in \fap \times \fa \times
\fa$ on the line segment joining $(X,T,S)$ and $(X', T', S')$ (so that
$(T'', S'')$ is well-situated and $X''$ lies in $\interior R_i(\delta',
T'',S'')),$ $X'' + \ker \Pi_1$ intersects a face $\tF$ of $R_i(\delta',
T'',S'')$ strictly contained in $F''$, the face corresponding for $F$. 
Notice also that since $F'' \supsetneq \tF$ satisfies (3.18),  there
must exist points $Y$ arbitrarily close to $X''$ such that $Y+\ker \Pi_1$
does not intersect $\tF$ (or even $\affspan \tF$).  Let us examine these
faces $\tF$. 

Let $\tF$ be a maximal face of $R_i(\delta', T,S)$ such that for some $X\in
\interior R_i(\delta', T,S)$ the intersection 
$(X+ \ker \Pi_1) \cap \affspan F$ is non-empty, but that
there exist points in any neighbourhood of $X$ such that the corresponding
intersection is empty.  We will call all such faces {\it problematic.} The
face $\tF$ must intersect $\ker \Pi_1$, so the projection $\oH$ of $\affspan
F$ to $\ofap$ is a subspace.  In fact, maximality of $\tF$
is easily seen to imply that $\oH$ is actually a hyperplane through $\oo$
in $\ofap$, and so can be given by an equality 
$$
\sum_{\lambda \in \cB} d_{\lambda}^{\tilde F} \lambda = 0. \eqno (3.19)
$$

An argument similar to that in Lemmas 3.4 and 3.5 shows that given a
problematic face $\tF$ of $\RDT$, and any well-situated $\Tp$, the
corresponding face $\tF'$ of $\RDTp$ is also problematic and its projection
to $\ofap$ is a hyperplane given by the equation (3.19). The
complement in $R_i(\delta', T,S)$ [resp. $\oRDT$] of the union  over all
problematic faces $\tF$ of $R_i(\delta', T,S)$, of the hyperplanes in $\fap$
[resp. $\ofap$] determined by (3.19), is a finite, disjoint union of convex
open polytopes, each of whose closures is given as the set of $X$ in
$R_i(\delta', T,S)$ [resp. $\oRDT$] satisfying 
$$
j(F) \sum_{\lambda \in \cB} d_{\lambda}^\tF \lambda (X) \geq 0, \quad
\text{for each problematic face } \tF, \eqno (3.20) 
$$
with $j$ an assignment of $\pm 1$ to
each problematic face of $R_i(\delta', T,S)$.  Let $J_i$ be the set of
assignments $j$ such that the set of $X$ in $\oRDT$ satisfying (3.20) has
non-empty interior in $\ofap$.  Then
$J_i$ is independent of well-situated $(T, S)$. 

\proclaim{Definition 3.5} (a) Given $(T,S)$ well-situated, $i\in I$,
$j\in J_i$, define $R_{i,j}(\delta',T,S)$ [resp.
$\oRjDT$] to be the set of $X$ in $R_{i}(\delta',T,S)$
[resp. $\oRDT$] satisfying  the inequalities (3.20).\hb
(b) Given $X\in\fa_P$, $i\in I$, define $R_i(T,S)_X$ to be the
polytope
$$
\bigl((X+\ker \Pi_1)\cap R_i(T,S)\bigr) -X\i \ker\Pi_1.
$$
\endproclaim

We have proven the following Lemma.

\proclaim{Lemma 3.6} Fix $i\in I$, $j\in J_i$.  Suppose that we are
given $(T,S),\ (T',S')$ well-situated, and 
$X\in R_{i,j}(\delta', T,S)$, $X'\in R_{i,j}(\delta', T',S')$.  If $Y$
is an extreme point of (3.14), let $F$ be a face of $R_i(T,S)$
satisfying (3.17), and let $F'$ be the corresponding face of
$R_i(T',S')$.  Then $(X'+\ker \Pi_1)\cap F'$ is an extreme point of
$(X'+\ker\Pi_1) \cap R_i(T',S')$.
\endproclaim

\proclaim{Corollary 3.7} Fix $i\in I$, $j\in J_i$, $\mu\in\fa^*$.
Then the integral 
$$
\int_{R_i(T,S)_X} e^{\mu(H)} dH, \quad (T,S) \text{ well-situated, }
X\in R_{i,j} (\delta',T,S),
$$
is t-finite in $X,T,S$.\endproclaim

\demo{Proof} This is immediate from Lemmas 3.6, 2.3, 2.4.\qed\enddemo

%%%%%%%Fix the next two paragraphs!!!
%%%%For each assignment $j \in J_i$, define $R_{i,j}(\delta', T,S)$
%%%%[resp. $\oRjDT$] to be the set of $X$ in $R_i(\delta', T,S)$
%%%%[resp. $\oRDT$] satisfying  the inequalities (3.20).  Then we
%%%%have shown that for any well-situated $(T,S)$ and $(T',S')$, any $X \in
%%%%\interior R_{i,j}(\delta', T,S)$ and $X' \in \interior \RjDT'$ and any face
%%%%$F$ of
%%%%$R_i(\delta', T,S)$ satisfying (3.18), (3.17) is an extreme point  of
%%%%(3.14) if and only if $(X' +\ker \Pi_1)\cap F'$ is an extreme point of
%%%%$(X'+\ker \Pi_1)\cap R_i(T', S')$, where $F'$ is the face of $R_i(\delta',
%%%%T',S')$ corresponding to $F$.  Since every face $F$ of $R_i(\delta', T,S)$
%%%%satisfying (3.18) is a face of $R_i(T,S)$ it can be given as the solution
%%%%set to a finite collection of equalities  from $(3.7)'$ and $(3.8)'$ and so
%%%%the extreme points of $(X+ \ker \Pi_1) \cap R_i(T,S)$ are linear in $X,T$
%%%%and $S,$ for $X\in \interior R_{i,j}(\delta', T,S)$. We state this as a lemma. 
%%%%
%%%%\proclaim{Lemma 3.5}Given $X\in \fap$ and $i \in I$, write $R_i(T,S)_X$ for
%%%%the polytope $((X+\ker \Pi_1) \cap R_i(T,S))-X \text{ in } \ker \Pi_1$.
%%%%For any choice of $j \in J_i$, the extreme points of $\RiTx$ are linear in
%%%%$X,T,$ and $S$, for $X\in R_{i,j}(\delta', T,S)$, and $\T$ well-situated.
%%%%\endproclaim

Notice lastly that for every $j,$ the polytope $\oRjDT$ is again a pyramid
with apex $\oo$, whose bounding half-spaces are given by the inequalities
from (3.15), (3.20) and (3.16) (giving the base), and so is independent of
$S$. 

Write $\overline{R_{i,j}}$ for the polyhedron bounded only by the
inequalities from (3.15) and (3.20); the penultimate step in the proof
of our main theorem (3.3) will be noticing that $\overline{R_{i,j}}$
does not depend on $T$ or $S$.

\beginsection 4. The Main Theorem.

In [8], we saw that if the rank of $G$ is at most two, then the truncated
integral $J_\fo^T(f,\pi)$ was asymptotic, as $T$ approached infinity in
certain sub-cones of $\fa^+$, to a t-finite function of $T$.  In this
section we prove this for general $G$.

\proclaim{Theorem 4.1} Let $G$ be a rational reductive group with
anisotropic centre, and let $\pi$ be a 
rational representation of $G$ on a finite-dimensional vector space $V.$
For each 
$\pi$-dependent cone $\C$ in $\fa^+$, each geometric equivalence class
$\fo \in \fO,$ and each Schwartz-Bruhat function $f$ on $V(\bA)$, there
exists a unique t-finite function $P_{\fo,\C}$ on $\fa$ such that for every
sufficiently small $\vare >0$ and every $c>0$
there exists a continuous seminorm $\|\cdot\|=\|\cdot\|_{\vare,c}$ on the
space of Schwartz-Bruhat functions on $V(\bA)$ such that
$$
\sum_{\fo\in\fO} \bigl| J_\fo^T(f,\pi) - P_{\fo,\C}(T) \bigr| <
\|f\|e^{-c\|T\|}, 
$$
for all $T$ in $\C_\vare.$ The
function 
$$ P_\C(T) = \sum_{\fo \in \fO} P_{\fo,\C}(T)$$
is also t-finite.
\endproclaim

The basic outline of our proof of this theorem is the same as that of
Theorem 6.1 of [8], however there are additional complications, arising from
the constructions of section 3 and the need to watch the dependence on $f$.
As in [8] we prove the Theorem by examining differences
$J_\fo(T+S)-J_\fo(T)$ and then applying Lemma 4.2, whose
proof is easily adapted from the proof of Lemma 4.2 of [8].

\proclaim{Lemma 4.2}  Suppose that $\{J_\fo\}_{\fo\in\fO}$ is a
collection of continuous functions on an open cone $\C$ in $\fa^+$, and
that there 
exists a collection, indexed by $\fo \in \fO,$ of t-finite functions
$\{p_\fo\}$ on $\fa\times\fa$ and a constant $b$
such that for every $c>0$ there exists a constant $C_c$ such that
$$
\sum_{\fo\in\fO} |J_\fo(T+S) - J_\fo(T) -p_{S,\fo}(T) | <
C_ce^{-c\|T\|},
$$
for all $S$ in $\C(1)$ and every $T$ in the cone with $\|T\|\ge b.$  Then
for each 
$\fo\in\fO$  there exists a unique t-finite function $P_\fo$ and for
every $c>0$ there exists a constant $d_c$ depending only on $c$ such that
$$
\sum_{\fo\in\fO} |J_\fo(T)- P_\fo(T) | < d_c C_c e^{-c\|T\|}
$$
for every $T$ in the cone.
\endproclaim

We used reduction theory to reduce the difference
$J_\fo(T+S)-J_\fo(T)$ to expressions of the form (1.3).  The idea now
is to use the Poisson summation formula and the constructions of
section 3 to reduce the problem to an application of Corollary 3.7.  A
simplified example to keep in mind is the following: for $f$ a
Schwartz function on $\bR$, the integral over $x$ in the
interval from $0$ to $T$ of $\sum_{n\in\bZ}f(e^x n)$ approaches
$$
(e^{-T}-1)\hf(0) + \int_0^\infty \sum_{n\ne 0} e^x \hf(e^x n)\,dx
$$
as $T$ tends to infinity in the cone $\bR^-$ of negative reals and approaches
$$
Tf(0) + \int_0^\infty \sum_{n\ne 0} f(e^x n)\,dx
$$
as $T$ tends to infinity in the cone $\bR^+$ of positive reals.  The
t-finite functions $T$ and $e^{-T}-1$ are integrals of the form of
Corollary 3.7. 

\demo{Proof of Theorem 4.1} Fix a geometric equivalence class $\fo \in \fO$
and an 
$\vare >0$ sufficiently small that the set $\C_\vare$ is
non-empty.  Clearly, if for each sufficiently small $\vare>0$ there exists
a t-finite function satisfying the above estimate, then it must be
unique and independent of $\vare.$  We need therefore only find an
approximation on the cone $\C_\vare$  for a fixed $\vare >0$.  Let us also
fix $c>0$.

To apply Lemma 4.2, we must
estimate the difference $J_\fo^{T+S}(f,\pi)- J_\fo^T(f,\pi).$
Notice that we can require $\|T\|$ to be larger than any constant that is
independent of $f.$  Throughout this proof this is what we will mean by
choosing $T$ sufficiently large. Take $T,S$ as in Lemma~4.2,
with $T$ sufficiently large. For the remainder of the proof we will
assume without mention that $(T,S)$ are well-situated, so that $T$ is
sufficiently large and lies in 
$\C_\vare$, and that $S$ lies in $\C_\vare(1).$  By a fixed constant
we mean one that is independent of $T$, $S$, $\fo$, and $f$.

Equation (1.3) states that
$$
\eqalignno{
J_\fo^{T+S}(f,\pi) - J_\fo^T(f,\pi) =& \sum_{P\sin Q\snot G}
\int_{N_P(\bQ)\\N_P(\bA)} \int_{A_P(\bR)^0} e^{-2\rho_P(H_P(a))}  
\phi_\fo (na,f^{P,K,T_2})\cr
&\qquad \times \Gamma_P^Q(H_P(a),T-T_2) \Gamma_Q(H_P(a) -T,S) dadn.\cr}
$$
Since the outer sum is finite, we need consider only the term corresponding
to a fixed choice of $P\sin Q\snot G.$
Since the function $f^{P,K,T_2}$ is Schwartz-Bruhat and is independent of
$T$ and $S$, we can reduce the problem to constructing a t-finite
function that approximates the integral
$$
\eqalignno{
&\int_{N_P(\bQ ) \\ N_P(\bA )}  \int_{A_P(\bR )^0}
e^{-2\rho_P (H_P(a))} \phi_\fo (na,f)\cr
&\qquad\quad \times\Gamma_P^Q (H_P(a), T- T_2 ) 
\Gamma_Q (H_P(a) -T,S) \,dadn,&(4.1)\cr}
$$
for a Schwartz-Bruhat function $f$ on $V(\bA ).$ 

Recall that the set $R_{P,Q}(T,S)$ defined in the previous section is the
support in $\fa_P$ of the characteristic function sending $X$ to
$$
\Gamma_P^Q (X , T ) \Gamma_Q (X-T,S).
$$
Therefore the integral over $a$ in (4.1) can be seen as an integral
over the set of $a$ such that $H_P(a)$ lies in $R_{P,Q}(T - (T_2)_P^Q,S).$
We spent a lot 
of effort in the previous section producing a decomposition of
$R_{P,Q}(T,S)$;
let us now make use of it.  Consider a closed region 
$R=R_i(T-(T_2)_P^Q,S)$ in the
decomposition 
$$
R_{P,Q}(T - (T_2)_P^Q,S) = \bigcup_{i\in I} R_i(T - (T_2)_P^Q,S)
$$
of $R_{P,Q}(T - (T_2)_P^Q,S)$ given in the previous section.   Since
$R_{P,Q}(T - (T_2)_P^Q,S)$ is 
the disjoint (modulo boundaries) union of these regions, and since
there are only finitely many of the regions, it is sufficient to
estimate the integral for $H_P(a)$ in this one region $R.$  Recall that
together with the region $R$ we associated a disjoint decomposition of the
set $\Pi$ of non-zero weights of $\pi$ with respect to $A_P$:
$$
\Pi = \Pi_- \cup \Pi_0 \cup \Pi_+;
$$
define the weight spaces
$$
V_- = \bigoplus_{\lambda\in\Pi_-} V^\lambda, \quad
V_0 = \bigoplus_{\lambda\in\Pi_0\cup\{0\}} V^\lambda, \quad
V_+ = \bigoplus_{\lambda\in\Pi_+} V^\lambda.
$$

Pick a basis of $V(\bQ)$ so that each basis element is in some
$V^\lambda$.  By setting the basis to be orthonormal,
we obtain an inner product $v\cdot w$ and a norm $\|v\|$ on $V(\bQ)$
and on $V(\bR)$---notice that for $v,w\in V(\bQ)$, $v\cdot w$ is
rational.  We also set $V(\bZ)$ to be the set of integral linear
combinations of the elements of this basis, and similarly with
$V({1\over N}\bZ)$ for any positive integer $N$.

Replace the integral over $N_P (\bQ)\backslash N_P(\bA )$ by one over
$\omega_P$, a relatively compact convex fundamental
domain of $N_P(\bQ )\\ N_P(\bA )$ containing the identity. With this
substitution, we need not worry about $N_P(\bQ)$-invariance of our
expressions. 

We must prove that there exists a t-finite approximation of the
integral
$$
\int_{\omega_P}  \int_{\exp R}  e^{-2\rho_P (H_P(a))} \phi_\fo
(na,f)\,dadn,
$$
where $\exp R=\{a\in A_P(\bR)^0\mid H_P(a)\in R\}$. 

Notice that because of Lemma 1.1, we have the following equality:
$$
\eqalignno{
\phi_\fo(na,f) =& \sum_{\gamma\in\fo} f(\pi(na)^{-1} \gamma)\cr
=&\sum_{\gamma \in (\fo \cap V_0 (\bQ )) + V_+(\bQ)} 
   f(\pi (na)^{-1} \gamma)  + \sum_{\gamma\in\fo\cap V(\bQ)'} 
   f(\pi(na)^{-1} \gamma), & (4.2)\cr}
$$
where $V(\bQ)'$ is the set of $\gamma\in V(\bQ)$ with a nonzero
component in $V_-.$

We claim that the second term of the right-hand side of (4.2) is an error
term.  More precisely, we claim that there exists a fixed continuous
seminorm
$\|\cdot\|_1$ on the space fixed of Schwartz-Bruhat functions on $V(\bA)$
such that the expression 
$$
\int_{\exp R} e^{-2\rho_P (H_P(a))} 
   \int_{\omega_P} \sum_{\gamma\in V(\bQ)'} 
   \bigl|f(\pi(na)^{-1} \gamma)\bigr| \,dnda, \eqno (4.3)
$$
can be bounded by  $\|f\|_1 e^{-c \| T \|}$ for $T$
sufficiently large.

We prove this as follows.  Define a function $\of$
on $V(\bA )$ by
$$
\of (v ) \ = \ \sup_{n\in \omega_P} \mid f(\pi (n^{-1} )v )|.
$$
Now, for 
$a \in A_P(\bR )^0 $ with $H_P (a)$ in 
$\fa^+_P$, the set $a^{-1} \omega_P a$ is contained in $\omega_P$, and
so, given any function $h$ on $N_P (\bA )$,
$$
\int_{\omega_P}  \ h(a^{-1} na )dn \ \leq \ \supn | h (n)|.  \eqno(4.4)
$$
(Notice that $\omega_P$ has volume one.)  Therefore, (4.3) is bounded by 
$$
\eqalignno{
&\int_{\exp R}  e^{-2\rho_P (H_P(a))} 
   \int_{\omega_P} \sum_{\gamma\in V(\bQ)'} 
   \bigl|f\bigl(\pi(a^{-1}na)^{-1} (\pi(a^{-1})\gamma)\bigr)\bigr|da\cr
\le & \int_{\exp R} e^{-2\rho_P (H_P(a))} \sum_{\gamma\in V(\bQ)'} 
    \of(\pi (a^{-1}) \gamma )da. &(4.5)\cr}
$$

Since $\omega_P$ is relatively compact, the function $\of$ is continuous
and rapidly
decreasing, that is, $\of$ is a finite sum of functions of the form
$\prod_v f_v$, with each $f_v$ a continuous function on $V(\bQ_v)$ that is
compactly supported if $v$ is finite, and decreases faster than the
inverse of any polynomial if $v$ is infinite.  In particular there is an
integer 
$$N_1(f)= \prod_p p^{n_p(f)}$$ 
determined by the support of $\of$ such that the sum in (4.5) can be taken
on 
$V({1\over N_1(f)} \bZ)' = V(\bQ)' \cap V({1\over N_1(f)}\bZ)$ instead of
$V(\bQ)'.$  

By the definition of $R$ we know that for $a$ in $\exp R$ and
$\lambda\in\Lambda_i^-$, we have 
$$
\lambda(H_P(a)) < -\delta_i B(T).
$$
Since $B$ is
positive on $\C$, $B$ is larger 
on all $\cC_{\vare}$ than some fixed multiple, depending on $\vare$, of the
norm. This implies the existence of a fixed positive constant $k$ 
depending on $\vare$ such that
$$
\lambda (H_P(a)) \ \leq \ -k \| T\| , \quad \hbox{for all } a\in \exp R,
\lambda\in \Pi_-.  \eqno (4.6)
$$
There is also a fixed positive constant $k'$, also depending on $\vare$, so
that
$$
\lambda (H_P(a)) \ \geq \ -k' \| T \|, \quad \hbox{for all } a\in \exp R,
\lambda\in \Pi_0 \cup\Pi_+.
$$
The inequality (4.6) implies that we can force all the points
$\pi(a^\1)\gamma,$ $a\in \exp R,\gamma\in V({1\over N_1(f)}\bZ)'$, to lie
outside any fixed compact set, by choosing $T$ sufficiently large.

Given a vector $v\in V(\bA)$ and a weight $\lambda\in \Pi$, write
$v^\lambda$ for the component of $v$ in $V^\lambda$
and $v_\bR$ for the component of $v$ in $V(\bR).$ Let $\mu$ be a weight in
$\Pi_-.$ 
Since $\of$ is rapidly decreasing the previous paragraph implies
that we can bound
$$
\eqalignno{
|\of (\pi(a^{-1}) \gamma ) |  &\leq \ 
   \threenorm f\threenorm_1 \| \pi (a^{-1}) \gamma_\bR^\mu \|^{-\ell} 
   \prod_{\lambda \in \Pi \setminus \mu} 
   (1+ \| \pi (a^{-1})\gamma_\bR^\lambda \|^n )^{-1}\cr
&=\threenorm f\threenorm_1 e^{\ell\mu(H_P(a))} \| \gamma_\bR^\mu\|^{-\ell}
   \prod_{\lambda\in \Pi\setminus \mu} (1 + e^{-n \lambda(H_P(a))}
\|\gamma_\bR^\lambda\|^n )^{-1} & (4.7)
\cr}$$
for all $\gamma\in V({1\over N_1(f)}\bZ)'$ with $\gamma^\mu$ nonzero, 
by choosing $T\in\C_\vare$ sufficiently large, 
where $\ell \ge n$ are arbitrarily large fixed integers.  The norm
$\threenorm\cdot\threenorm_1$ is given by
$$
\threenorm f\threenorm_1=\sup_{\mu\in\Pi_- \atop v\in V(\bA), v_\bR^\mu\ne
0}
   \Bigl({\of(v)    \| v_\bR^\mu \|^{\ell} 
   \prod_{\lambda \in \Pi \setminus \mu} 
   (1+ \| v_\bR^\lambda \|^n )}\Bigr),
$$
and is continuous with respect to the topology on the space of
Schwartz-Bruhat functions on $V(\bA).$

For $\lambda\in \Pi_-$ and $x$ a non-negative real number we can bound 
$$
(1 + e^{-n \lambda(H_P(a))} x )^{-1} \le (1 + x )^{-1}.
$$
For $\lambda \in \Pi_0\cup\Pi_+$ and $x$ a non-negative real number we can
bound 
$$
\eqalign{
1 + e^{-n \lambda(H_P(a))} x  \ge& \min (1, e^{n\lambda(H_P(a))}) (1 +
x)\cr
\ge & e^{-nk' \|T\|} (1 + x).\cr}
$$
Also, since $\mu$ lies in $\Pi_-,$
$$
e^{\ell \mu(H_P(a))} \le e^{-\ell k \|T\|}.
$$
Putting together the above inequalities, we obtain that $\of(\pi(a^\1)\g)$
is bounded by 
$$
\threenorm f\threenorm_1\exp\bigl( - (\ell k - |\Pi_0\cup\Pi_+| nk')
   \|T\|\bigr)
   \|\gamma_\bR^\mu\|^{-\ell} \prod_{\lambda \in\Pi\setminus \mu} 
   (1 + \|\gamma_\bR^\lambda\|^n  )^{-1},
$$
where $|\Pi_0\cup\Pi++$ denotes the cardinality of the set.  Since
$\|\g_\bR^\mu\|\ge 1/N_1(f)$ with $N_1(f)$ a positive integer, and
$\ell\ge n$, we have
$$
\|\gamma_\bR^\mu\|^{-\ell}=
\|\gamma_\bR^\mu\|^{-(\ell-n)}\|\gamma_\bR^\mu\|^{-n}
\le N_1(f)^{\ell-n}{(N_1(f)^n +1)\over 1 + \|\g_\bR^\mu\|^n}
\le 2 N_1(f)^\ell (1 + \|\g_\bR^\mu\|^n)^\1,
$$
and so we obtain the bound
$$
2 N_1(f)^\ell \threenorm f\threenorm_1 
\exp\bigl( - (\ell k - |\Pi_0\cup\Pi_+|nk') \|T\|\bigr) 
   \prod_{\lambda\in\Pi} (1 + \|\gamma_\bR^\lambda\|^n  )^{-1}
$$
of $\of(\pi(a^\1)\g)$; this bound is independent of the weight $\mu$,
and hence is valid for all $\gamma\in V({1\over N_1(f)}\bZ)'.$  
Choose $n$ so large that for every $\lambda \in \Pi$, the sum 
$$
\sum_{\gamma \in V^\lambda ( \bZ )} (1 + \| \gamma_\bR \|^n )^{-1}
$$
converges, and write $C_\lambda$ for its value.  Then for every natural
number $N$ and every $\lambda \in \Pi$, the sum  
$$
\sum_{\gamma \in V^\lambda ({1\over N} \bZ )} (1 + \|
\gamma_\bR \|^n )^{-1}
$$
is bounded by $C_\lambda N^n.$
Choose $\ell\ge n$
so that $\ell k - | \Pi_0 \cup \Pi_+ | nk' \ \ge \ c+1.$ 
Then (4.5) is bounded by
$$
\eqalign{
& 2 N_1(f)^\ell \threenorm f\threenorm_1  e^{-(c+1)\| T \|} 
\Bigl(\int_{\exp R}  e^{-2\rho_P (H_P(a))}da \Bigr)  
\sum_{\gamma \in
V({1\over N_1(f)} \bZ )} \prod_{\lambda \in \Pi} (1 + \| \gamma_\bR^\lambda
\|^n )^{-1}\cr
& \quad \leq \ 2 N_1(f)^\ell \threenorm f\threenorm_1 e^{-(c+1)\| T \|}
\volume (R)
\prod_{\lambda \in \Pi} 
\sum_{\gamma \in V^\lambda ({1\over N_1(f)} \bZ )} (1 + \|
\gamma_\bR\|^n)^{-1}
\cr
&\quad \le C_1 N_1(f)^{\ell+|\Pi|n} \threenorm f\threenorm_1  e^{-(c+1)\| T
\|} 
\volume (R),\cr}
$$
where 
$$
C_1=2 \prod_{\lambda\in\Pi} C_\lambda
$$
is a fixed constant. By Lemma 3.4, 
the extreme points of the region $R$ are linear in $T - (T_2)_P^Q$ and $S$
for all  
sufficiently large $T\in \C_\vare$ and all $S\in \C_\vare(1)$, so the
volume of $R$ is a polynomial in  $T$ and $S.$  Therefore 
$$
C_1 N_1(f)^{\ell+|\Pi|n} \threenorm f\threenorm_1 e^{-(c+1)\| T \|} \volume
(R) \
\leq \ C_1' N_1(f)^{\ell+|\Pi|n} \threenorm f\threenorm_1 e^{-c \|T\|} 
$$
for some fixed constant $C_1'.$  This completes our bound of (4.3), since
the norm $\|\cdot\|_1$ given by $\|f\|_1 = C_1' N_1(f)^{\ell+|\Pi|n} 
\threenorm f\threenorm_1$ is a fixed continuous seminorm.

Rewrite the first term of (4.2), with our fixed geometric equivalence class
$\fo$, as 
$$
\eqalignno{
& \sum_{\gamma_0 \in \fo \cap V_0 (\bQ )}\ \sum_{\gamma_+ \in V_+ (\bQ )}
f(\pi (na)^{-1} (\gamma_0 + \gamma_+ ))\cr
\ = & \sum_{\gamma_0 \in \fo \cap V_0 (\bQ )}\ \sum_{\gamma_+ \in V_+ (\bQ
)}
\int_{V_+ (\bA )}  f(\pi (na)^{-1} (\gamma_0 + v ))\cdot 
\psi (\gamma_+ \cdot v)dv \cr
\ = & \sum_{\gamma_0 \in \fo \cap V_0 (\bQ )} 
\int_{V_+ (\bA )}  f(\pi (na)^{-1} (\gamma_0 + v ))dv & (4.8)\cr
& \quad +  \sum_{\gamma_0 \in \fo \cap V_0 (\bQ )}\ \sum_{\gamma_+ \in V_+
(\bQ )\setminus \{0\}} 
\int_{V_+ (\bA )} f(\pi (na)^{-1} (\gamma_0 + v )) 
\psi (\gamma_+ \cdot v)dv ;\cr}
$$
the first step was Poisson summation on $V_+ (\bQ ).$  We have chosen
here $\psi$ to be the standard additive character on $\bA$ given in
Tate's thesis.

The following expression dominates the sum over all geometric equivalence
classes $\fo$ of the second term of the right-hand side of (4.8):
$$
\sum_{\gamma_0 \in  V_0 (\bQ )}\ 
\sum_{\gamma_+ \in V_+ (\bQ )\setminus \{0\}}
\Bigl|\int_{V_+ (\bA )}  f(\pi (na)^{-1} (\gamma_0 + v )) 
\psi (\gamma_+ \cdot v)dv\Bigr|
\eqno (4.9) $$
We claim that the integral over $a$ and $n$  of (4.9) is another error
term.  We prove this as follows.  

Define the function
$$
f_{a,n} (\gamma_0 , \gamma_+ ) \ = \ \int_{V_+ (\bA )} f(\pi (an)^{-1}
(\gamma_0 + v)) \psi (\gamma_+ \cdot v)dv ,\quad \gamma_0 \in V_0 (\bA ),
\ \gamma_+ \in V_+ (\bA ).
$$
For $\gamma_0 \in  V(\bQ), \gamma_+\in V_+\setminus \{0\},$
$|f_{a,n}(\g_0 , \g_+)|$ 
is the summand in (4.9) corresponding to $\gamma_0$ and $\gamma_+.$
By (4.4), we have the inequality
$$
\int_{\omega_P} |f_{a,n} (\gamma_0,\gamma_+)| dn 
\le \sup_{n\in \omega_P} \ | \ f_{a,n} (\gamma_0 , \gamma_+)\ |.
$$
Notice that
$$
\eqalign{
f_{a,n} (\gamma_0,\gamma_+) \ = & \ \int_{V_+ (\bA )} 
f(\pi (n^{-1}) \pi (a^{-1} )(\gamma_0 + v)) \psi (\gamma_+ \cdot v)dv\cr
\ = & \ \int_{V_+ (\bA )}  f(\pi (n^{-1})(\pi (a^{-1} )
\gamma_0 +\pi (a^{-1}) v)\psi (\pi (a)\gamma_+ \cdot\pi(a^{-1}) v)dv\cr
= & e^{2\rho_+ (H_P(a))} \int_{V_+ (\bA )} f(\pi (n^{-1}) \pi (a^{-1} )
\gamma_0 + v)) \psi (\pi(a)\gamma_+ \cdot v)dv \cr
= & e^{2\rho_+ (H_P(a))} \ f_{1,n}(\pi (a^{-1})\gamma_0,\pi
(a)\gamma_+),\cr}
$$
where $e^{2\rho_+(H_P(a))}$ is the Jacobian of the change of variables
$$
v \ \mapsto \ \pi (a) v, \quad v\in V_+(\bA)
$$
on $V_+ (\bA )$, for $a \in A_P (\bR )$.  Now,
the function $f_{1,n}$ is  
Schwartz-Bruhat for each $n\in \omega_P$ and continuous with respect to $n.$
Since $\omega_P$ is relatively compact the function $f_N$ on $V_0(\bA)
\times V_+(\bA)$ defined by
$$
f_N (v_0 ,v_+ ) \ = \ \sup_{n\in \omega_P} \ | f_{1,n} (v_0 , v_+ )|,
\quad v_0\in V_0(\bA), v_+\in V_+(\bA)
$$
is continuous and rapidly decreasing.

The integral over $a$ and $n$ of (4.9) is therefore bounded by
$$
\int_{\exp R} \ e^{(2\rho_+ - 2\rho_P) (H_P(a))} \sum_{\gamma_0 \in V_0 (\bQ
)}
\sum_{\gamma_+ \in V_+ (\bQ )\setminus \{0\}} f_N (\pi (a^{-1} )
\gamma_0 , \pi (a) \gamma_+)da, \eqno (4.10)
$$
where, since the integrand is positive and the expression clearly
converges, Fubini's theorem allowed the free interchange of integrals.
Fubini's theorem will trivially apply through (4.23) because all integrals
will be over compact sets of continuous functions and its use will not be
mentioned.  The function $f_N$ is rapidly decreasing, $a\in \exp R$ is real,
so the sums in  (4.10) can be taken in $V({1\over N_2(f)} \bZ) $ instead of
$V(\bQ),$ for some integer $N_2(f)$ depending only on $f.$

Pick new fixed positive constants $k,k'$ (depending on $\vare$) 
so that
$$
\eqalign{
\lambda (H_P(a)) \ \geq \ k \| T \|,\qquad \text{for every } a \in \exp R,\
\lambda \in \Pi_+, \cr
\lambda (H_P(a)) \ \leq \ k' \| T \|,\qquad \text{for every } a \in \exp R,\
\lambda \in \Pi_0; \cr}
$$
this again is possible because of the inequalities in (3.7).
Let $a\in \exp R$ and $\g_+\in V_+({1\over N_2(f)} \bZ) \setminus \{0\},$ be
arbitrary, and let $\mu\in\Pi^+$ be any weight in $\Pi_+$ such that
$\g_+^\mu$ is non-zero.  Then
$$
\| (\pi (a)\gamma_+ )_\bR \| \ \geq \ \| \pi (a) (\gamma_+^\mu )_\bR
\| \ = \ e^{\mu (H_P(a))} \| (\gamma_+^\mu )_\bR \| \ \geq \
{1\over N_2(f)} e^{k \|T\|}.
$$
Since $f_N$ is rapidly decreasing, the following inequality holds for
arbitrarily large integers  $\ell \geq n$, $a\in \exp R,\;\g_0\in
V_0({1\over N_2(f)}\bZ)$, and  $\g_+\in V_+({1\over N_2(f)} \bZ) \setminus
\{0\}$, if
$T\in\C_\vare$ is sufficiently large:
$$
\eqalign{
&| f_N (\pi (a^{-1} ) \gamma_0 , \ \pi (a) \gamma_+ ) | \cr 
\leq \ & \threenorm f\threenorm_2 \, \| \pi (a) (\gamma_+ )_\bR^\mu
\|^{-\ell} 
 \prod_{\lambda \in \Pi_0}
 (1 + \| \pi (a^{-1} ) (\gamma_0 )_\bR^\lambda \|^n )^{-1}
 \prod_{\lambda \in \Pi_+ \backslash \mu} (1 + \|
 \pi (a) (\gamma_+ )_\bR^\lambda \|^n )^{-1}\cr
\leq \ & 2 N_2(f)^\ell
 \threenorm f\threenorm_2 \exp \bigl(-\|T\| (\ell k - |\Pi_0 | nk')\bigr) 
 \prod_{\lambda \in \Pi_0 \cup \Pi_+} 
(1 + \| (\gamma_0 + \gamma_+ )_\bR^\lambda \|^n )^{-1},\cr}
$$
where $\threenorm \cdot\threenorm_2 $ is a fixed continuous seminorm.
Write $\ell' = \ell k - |\Pi_0|nk'$.  We have bounded the
summand in (4.10) independently of the choice of $\mu$, so (4.10) is
bounded by
$$
\eqalign{
&  2N_2(f)^{\ell} \threenorm f \threenorm_2 e^{-\ell' \|T\|} 
   \left( \int_{\exp R} e^{(2\rho_+-2\rho_P)(H_P(a))}da\right) 
   \sum_{\g \in (V_0 \oplus V_+)(\frac{1}{N_2(f)} \bZ)} 
   \prod_{\lambda \in \Pi_0 \cup \Pi_+}
   (1+\|\g_{\bR}^{\lambda}\|^n)^{-1}\cr
&= 2N_2(f)^{\ell} \threenorm f \threenorm_2 e^{-\ell' \|T\|} 
   \left( \int_{\exp R} e^{2\rho_+-2\rho_P)(H_P(a))}da\right) 
   \prod_{\lambda \in \Pi_0 \cup \Pi_+} \sum_{\g^{\lambda}\in V^{\lambda} 
   (\frac{1}{N_2(f)} \bZ)} (1+\|\g_{\bR}^{\lambda}\|^n)^{-1}.\cr}
$$
Choose $n$ so large that the sum 
$$ 
\sum_{\gamma^{\lambda} \in V^\lambda (\bZ )} (1 + \| \gamma^{\lambda}_\bR
\|^n
)^{-1}
$$
converges for every $\lambda \in \Pi_0 \cup \Pi_+$, so that (4.10) is
bounded by
$$
\eqalignno{
& \quad  \le \ C_2 N_2(f)^{\ell+|\Pi_0\cup\Pi_+|n} 
 \threenorm f\threenorm_2 e^{-\ell' \|T\|} 
\Bigl(\int_{\exp R}  e^{(2\rho_+ - 2\rho_P)(H_P(a))}da \Bigr),  &(4.11)}
$$
for some fixed constant $C_2.$  Since the extreme points of $R$ are linear
in $T-(T_2)_P^Q$ and $S$,  
the expression in parentheses in (4.11) is bounded by the exponential of
some fixed multiple of $\|T\|.$  By choosing $\ell$ sufficiently large we
obtain that (4.10) is bounded by $C_2' N_2(f)^{\ell+|\Pi_0\cup\Pi_+|n}
 \threenorm f\threenorm_2 e^{-c \|T\|}$ for all sufficiently
large $T$, where $C_2'$ is some fixed constant.  The seminorm
$\|\cdot\|_2$ given by $\|f\|_2 = C_2' N_2(f)^{\ell+|\Pi_0\cup\Pi_+|n}
 \threenorm f\threenorm_2$ is continuous, so that this term, too, is
an error term. 

We can now deal with the integral of the first summand
of (4.8),
$$
\int_{\exp R}  e^{-2\rho_P (H_P(a))} \int_{\omega_P} \sum_{\gamma
\in \fo \cap V_0 (\bQ )} \int_{V_+ (\bA )} 
f(\pi (na)^{-1} (\gamma + v))dvdnda. \eqno (4.12)
$$
We are trying to prove that this integral has a t-finite approximation.

Given a subset $\cS$ of $\Pi_0$, write $\fS(\cS)$ for the set of weights in
$\Pi$ that vanish on all 
$$
(\ker \cS) \cap R = \ker (\sum_{\lambda \in \cS}(\sgn \lambda)\lambda)\cap
R;
$$
by (3.7) we have $\fS(\cS) \sin\Pi_0$.

Given $v\in V$, write
$\supp v$ for the set of those weights $\lambda$ such that $v$ has
a non-zero component in the weight space $V^\lambda.$  Given a subset
$\Pi'$ of $\Pi_0$, write
$$
\displaylines{
W_0(\Pi') = \bigoplus_{\lambda\in\Pi'} V^\lambda, \qquad
W_+(\Pi') = \bigoplus_{\lambda\in\Pi_0^+\setminus\Pi'} V^\lambda, \cr
W_0'(\Pi') = \{ v\in W_0(\Pi')
\mid \fS(\supp(v) \cap \Pi_0^-) = \Pi' \}. \cr}
$$
Notice that $W_0'(\Pi') $ is empty if $\fS(\Pi')\ne \Pi'$.
By Lemma 1.1, we know that
$$
\fo\cap\bigl(W_0'(\Pi')(\bQ) + W_+(\Pi')(\bQ)\bigr) =
\bigl( \fo\cap W_0'(\Pi')(\bQ)) + W_+(\Pi')'(\bQ).
$$
We can write
$$
\sum_{\g\in\fo\cap V_0(\bQ)} F(\g) = 
\sum_{\Pi'\sin \Pi_0\atop \fS(\Pi')=\Pi'}
\sum_{\g\in\fo\cap W_0'(\Pi')(\bQ)} 
\sum_{\g_+ \in W_+(\Pi')(\bQ)} F(\g + \g_+),
$$
for any function $F$ on $V_0(\bQ)$; a vector $\g \in \sigma \cap V _0(\bQ)$
appears in the summand corresponding to $\Pi' = \fS (\supp \g \cap \Pi_0^-).$

There are only finitely many choices
for $\Pi'.$  Make one and write $W_0=W_0(\Pi')$, $W_+=W_+(\Pi')$, and
$W_0'=W_0'(\Pi').$  Consider the summand in
(4.12) corresponding to $\Pi'.$  It is
$$
\eqalignno{
&\int_{\exp R} e^{-2\rho_P(H_P(a))} \int_{\omega_P} \sum_{\g\in\fo\cap
W_0'(\bQ)}
\sum_{\g_+\in W_+(\bQ)} 
\int_{V_+(\bA)} f\bigl( \pi(na)^\1 (\g + \g_+ + v) \bigr)
\,dv\,dn\,da\cr
=& \int_{\exp R} e^{-2\rho_P(H_P(a))} \int_{\omega_P} \sum_{\g\in\fo\cap
W_0'(\bQ)}
\sum_{\g_+\in W_+(\bQ)} & (4.13) \cr
&\qquad \int_{W_+(\bA) \oplus V_+(\bA)} 
f\bigl( \pi(na)^\1 (\g + v) \bigr) \psi(\g_+ \cdot v)\,dv\,dn\,da.\cr}
$$
We will next break up the sum over $\g_+$ in (4.13).

Given a subset $\cS$ of $\Pi_0^+$, write $\fb(\cS)$ for the set of
weights in $\Pi_0^+$ that vanish on all 
$$
\ker(\Pi'\cup \cS)\cap R=\ker (\sum_{\lambda \in \Pi' \cup \cS} 
(\sgn \lambda )\lambda )\cap R.
$$  
Given a subset $\Pi''$ of $\Pi_0^+$, write
$$
\displaylines{
U_0(\Pi'')= \bigoplus_{\lambda\in \Pi''\setminus \Pi'} V^\lambda, \quad
U_+(\Pi'')= V_+ \oplus \bigoplus_{\lambda\in \Pi_0^+\setminus\Pi''}
V^\lambda, \cr 
U_0'(\Pi'') = \{v\in U_0(\Pi'') \mid \fb(\supp v) = \Pi''\}.\cr}
$$
We can write
$$
\displaylines{
\sum_{\g_+ \in W_+(\bQ)} F(\g_+)=
\sum_{{\Pi'' \subset \Pi_0^+}\atop{\fb (\Pi'')=\Pi''}} 
\sum_{\g_+ \in U'_0 (\Pi'')(\bQ)} F(\g_+),}
$$
for any function $F$ on $W_+(\bQ)$; a vector $\g_+ \in W_+ (\bQ)$ appears
in the summand corresponding to $\Pi''= \fb (\supp \g_+)$.  

There are only
finitely many choices for $\Pi''$.  Make one and write $U_0=U_0(\Pi'')$,
$U_+=U_+(\Pi'')$, and $U_0' = U_0'(\Pi'').$  The summand of (4.13)
corresponding to it is 
$$
\eqalignno{
&\int_{\exp R} e^{-2\rho_P(H(a))} \int_{\omega_P} \sum_{\g\in\fo\cap
W_0'(\bQ)} \sum_{\g_+\in U_0'(\bQ)} & (4.14) \cr
&\quad \int_{W_+(\bA) \oplus V_+(\bA)} 
f\bigl( \pi(na)^\1 (\g + v) \bigr) \psi(\g_+ \cdot v)\,dv\,dn\,da.\cr}
$$

Let $\Pi_1\sin \Pi_0$ be the set of weights in $\Pi$ that vanish on all
$\ker(\Pi' \cup \Pi'') \cap R.$ This latter polytope equals 
$$
\ker\Bigl(\sum_{\lambda\in\supp(\g+\g_+)}(\sgn\lambda)\lambda\Bigr) \cap R. 
\eqno (4.15)
$$
for any $\g \in W_0', \g_+ \in U_0'$. Then 
$$
\ker \Pi_1 \cap R = \ker(\Pi'\cup\Pi'') \cap R,
$$
so $\Pi_1$ also equals the set of weights in $\Pi_0$ that vanish on
$\ker\Pi_1\cap R.$  Recall that in the previous section we defined a number
$\delta'$ and a decomposition
$$
R_i(\delta',T-(T_2)_P^Q,S) = \bigcup_{j\in J} R_{i,j}
(\delta',T-(T_2)_P^Q,S)
$$
of part of the region $R=R_i(T-(T_2)_P^Q,S)$, that depended on $\Pi_1$ so
that the regions $\overline{R_X} = \overline{R_i(T,S)_X}$ defined in Lemma
3.5 behaved well for $X$ in any given $R_{i,j}(\delta', T,S)$.  We will
soon use this decomposition.

First, we break up the integral over $a$ as follows. Let $\oA_P (\bR)^0
\subset A_P (\bR)^0$ be a complement to the subgroup $\exp (\ker \Pi_1)$ so
that the natural projection $p: \oA_P (\bR)^0 \mapsto \ofap = \fap/\ker
\Pi_1$ is an isomorphism, and normalize Haar measures on these two
subgroups so that their product is $da$.  We do this
independently of $f, T, S$.  Write $a=a_0 \oa,\ a_0 \in \exp (\ker
\Pi_1), \oa\in \oA_P (\bR)^0,$ for the canonical decomposition and write
$\exp\oR = \{\oa \mid H_P(a) \in R\}$.  Then the integral (4.14) becomes 
$$
\eqalignno{
&\int_{\exp\oR} e^{-2\rho_P(H_P(\oa))} \int_{\exp R_{H_P(\oa)}} 
  e^{-2\rho_P(H_P(a_0))} \int_{\omega_P} \sum_{\g\in\fo\cap W_0'(\bQ)} 
  \sum_{\g_+\in U_0'(\bQ)} & (4.16)\cr
&\quad \int_{W_+(\bA) \oplus V_+(\bA)} 
f\bigl( \pi(na_0 \oa)^\1 (\g + v) \bigr) \psi(\g_+ \cdot v)
\,dv\,dn\,da_0\, d\oa.}
$$
where $\exp R_{H_P(\oa)}$ is the exponential of the set $R_{H_P(\oa)} = R_i
(T-(T_2)_P^Q, S)_{H(\oa)}$ defined in Lemma 3.5. We will see that the
dependence of the innermost integral on $a_0$ is particularly simple. 

We first show that $N_P$ preserves both $U_+ (\bA)$ and $W_0 \oplus W_+
\oplus V_+$.  The two facts  are proven similarly, so we will consider only
the second.  It is clearly sufficient to prove that $\pi (n) V^{\lambda}$
lies in $W_0 \oplus W_+ \oplus V_+ = \oplus_{\mu \in \Pi' \cup \Pi_0^+ \cup
\Pi_+} V^{\mu}$ for every vector $v^{\lambda} \in V^{\lambda}, \lambda \in
\Pi' \cup \Pi_0^+ \cup \Pi_+$.  A weight $\mu$ in 
$\supp (\pi (n)v^{\lambda})$ is the sum of $\lambda$ and a non-negative
linear
combination of $\alpha \in \Delta_P$, and so by the inequalities (3.7),
must be at least as large as $\lambda$ on $R$.  If $\mu$ is 0 or 
$\sgn \mu =1$ we have nothing to show, while if $\sgn \mu = -1$, then  for
any $X\in (\ker \Pi') \cap R, \mu (X)$ must be both at most (since 
$\mu \in \Pi^-$) and at least (since $\mu (X) \geq \lambda (X)$) zero, so
that $\mu$ lies in $\Pi'$, proving the fact. 

We now prove that for $a_0 \in \exp(\ker \Pi_1), n \in N_P (\bA),$ and $w
\in W_0(\bA) \oplus W_+ (\bA) \oplus V_+ (\bA) = W_0 (\bA) \oplus U_0 (\bA)
\oplus U_+ (\bA),$ 
$$
\pi(na_0 n^\1) w -w \hbox{ lies in } U_+(\bA). \eqno (4.17)
$$
We have already shown that $\pi(n^\1)w$ can be written in the form 
$\pi(n^\1)w =w_0+ w_+$, with $w_0\in W_0(\bA) + U_0(\bA), w_+\in
U_+(\bA).$  The action of $a_0$ on $W_0 \oplus U_0$ is trivial, so
$$
\pi(a_0n^\1) w = \pi(a_0)(w_0 + w_+) = w_0 + \pi(a_0)w_+,
$$
and hence
$$
\pi(a_0n^\1) w - \pi(n^\1) w = \pi(a_0) w_+ - w_+ \in U_+(\bA).
$$
Since the action of $N_P$ preserves $U_+(\bA)$, we have proven (4.17) which
implies that the change of variables 
$$
v \mapsto \pi(na_0n^\1) v + \bigl(\pi(na_0n^\1) \g - \g\bigr)
$$
is a isomorphism on $W_+(\bA)\oplus V_+(\bA)$ that does not change
$\g_+\cdot v$ for any $\g_+\in U_0'(\bQ)$; its Jacobian is
$e^{2\rho_+(H_P(a_0))}$, 
where $2\rho_+$ is the sum of all weights in $\Pi^+\setminus\Pi'$,
including multiplicities.  The integral (4.16) therefore equals 
$$
\eqalignno{
&\int_{\exp\oR} e^{-2\rho_P(H_P({\oa}))} 
\Bigl(\int_{R_{H_P(\oa)}} e^{(2\rho_+' - 2\rho_P)(Y)} dY\Bigr)
\int_{\omega_P} \sum_{\g\in\fo\cap W_0'(\bQ)}\cr
&\quad \sum_{\g_+\in U_0'(\bQ)} \int_{W_+(\bA) \oplus V_+(\bA)} 
f\bigl( \pi(n\oa)^\1 (\g + v) \bigr) \psi(\g_+ \cdot v)\,dv\,dn\,d\oa,\cr}
$$

The region $R=R_i(T-(T_2)_P^Q,S)$ breaks up as 
$$
R=\bigcup_{j\in J_i} R_{i,j}(\delta', T-(T_2)_P^Q,S) \cup (R \setminus
R_i(\delta', T-(T_2)_P^Q,S)), 
$$
where the sets on the right are disjoint modulo boundary; this gives a
similar decomposition of $\exp\oR$.  We must estimate the contribution to
(4.16) of the integral over $\oa$ in each piece of the decomposition of
$\exp\oR$. 

We first claim that the last piece in this decomposition gives an error
term, that is, that the integral over $\oa$ and $Y$ with 
$H_P(\oa) + Y\in (R \setminus R_i(\delta', T-(T_2)_P^Q,S))$
of
$e^{-2\rho_P(H_P({\oa}))} e^{(2\rho_+' - 2\rho_P)(Y)}$ times  
$$
\int_{\omega_P} \sum_{\g\in\fo\cap W_0'(\bQ)}
\sum_{\g_+\in U_0'(\bQ)} \left| 
\int_{W_+(\bA) \oplus V_+(\bA)} 
f\bigl( \pi(n\oa)^\1 (\g + v) \bigr) \psi(\g_+ \cdot v)\,dv\right| dn
\eqno (4.18) 
$$
can be bounded by $\|f\|_3 e^{-c \|T\|}$ for some fixed continuous seminorm
$ \| \ \|_3$.  As with (4.3) and (4.9), we can bound (4.18) by an
expression of the form 
$$
e^{2\rho_P'(H_P({\oa}))} \sum_{\g \in W_0' ({1\over N_3(f)} \bZ )}
\sum_{\g_+ \in U_0' ({1\over N_3(f)} \bZ)} f_N(\pi (\oa)^{-1} \g, \pi (\oa)
\g_+), 
$$
with $f_N$ a continuous, rapidly decreasing function on $W_0 (\bA) \times
U_0 (\bA),$ and $N_3(f)$ an integer determined by the support of $f$.
Also, there exists a fixed constant $k''$ such that
$$
e^{(2\rho_+'- 2\rho_P)(H_P({\oa}))}\int_{R_P(\oa)} e^{(2\rho_+' -
2\rho_P)Y} dY
$$
is bounded by $e^{k'' \|T\|}$ for any $\oa \in \exp\oR$, for all
sufficiently large $T\in \cC_{\vare}$ (recall that $R_{H_P (\oa)}$ is a
slice in $R=R_i(T-(T_2)_P^Q,S)$).

Points $X$ in $R \setminus R_i(\delta', T-(T_2)_P^Q,S)$ satisfy
$$
\sum_{\lambda \in \cB} (\sgn \lambda) \lambda (X) \geq \delta' B(T-(T_2)_P^Q),
\eqno (4.19)
$$
with $\cB \subset \Pi_1$ a previously selected basis of $\spn \Pi_1$.
Since $\Pi_0$ is finite and $T$ is sufficiently large,
there exists a new fixed (and hence independent of $\Pi_1$) constant $k$,
depending on $\vare$, so that for each such $X$, some $\mu \in \Pi_1$
satisfies $(\sgn \mu) \mu (X) \geq k \|T\|$.  

Now, let $X\in R \setminus R_i(\delta', T-(T_2)_P^Q,S)$, $\mu$ as in
the previous paragraph, and $\g$ be any vector in $W_0'(\bQ), \g_+$ in
$U_0' (\bQ)$.  Since $\mu$ lies in $\Pi_1$, it vanishes on (4.15).
The boundary hyperplanes of $R$ that contain (4.15) are all of the form
$(3.7)'$ or $(3.8)'$; let $\lambda_1, \ldots, \lambda_k$ be corresponding
functionals (so that each $\lambda_i$ lies in $\Pi \cup \Delta_P \cup
\Delta_Q \cup \hD_Q \cup \hD_P^Q$).  The theory of the polar (see [5],
Theorem 6.4) implies that
$$
\sum_{\lambda \in \supp (\g +\g_+)} (\sgn \lambda) \lambda 
\eqno (4.20)
$$
is a positive linear composition of $(\sgn \lambda_i) \lambda_i$, and that
$(\sgn \mu) \mu$ is a non-negative linear combination of $(\sgn \lambda_i)
\lambda_i$.  All the constants involved in these linear combinations can be
chosen independently of $T$ and $S$, since all the functionals lie in the
finite set $\Pi \cup \Delta_P \cup \Delta_Q \cup \hD_Q \cup \hD_P^Q$.
Therefore, $(\sgn \mu) \mu$ is at most a fixed multiple of (4.20) and
so there exists  a new fixed constant $k'$, depending on $\vare$, such that
some $\lambda' \in \supp (\g+\g_+)$ satisfies 
$$
(\sgn \lambda') \lambda' (X) \geq k' \|T\|.
$$

At this point, we continue as with (4.10) to prove our claimed bound of
(4.18).  Therefore we need only estimate the contribution to (4.16) of
the integral over $\oa$ of the regions 
$$
\exp \Bigl(\,\overline{R_{i,j}(\delta', T-(T_2)_P^Q,S)}\,\Bigr)
 = \{ \oa \mid H_P(a) \in R_{i,j}(\delta', T-(T_2)_P^Q,S) \} \eqno (4.21)
$$
for each $j\in J_i$.

Fix $j\in J_i$, and write $\exp \overline{R}_j(T)$ for the set (4.21),
where we include the $T$ to remind ourselves of the dependence on $T$.  The
contribution of $\exp\overline{R}_j(T)$ is
$$
\eqalignno{
&\int_{\exp\oR_j(T)} e^{-2\rho_P(H_P({\oa}))} 
\Bigl(\int_{R_{H_P(\oa)}} e^{(2\rho_+ - 2\rho_P)(Y)} dY\Bigr)
\int_{\omega_P} \sum_{\g\in\fo\cap W_0'(\bQ)} &(4.22)\cr
&\quad \sum_{\g_+\in U_0'(\bQ)} \int_{W_+(\bA) \oplus V_+(\bA)} 
f\bigl( \pi(n\oa)^\1 (\g + v) \bigr) \psi(\g_+ \cdot v)\,dv\,dn\,d\oa,\cr}
$$

Corollary 3.7 says that the function
Lemma 3.5 says that for $\oa \in \exp\oR_j(T)$, the extreme points of
$R_{H_P(\oa)}$ are linear in $H_P(\oa), T-(T_2)_P^Q$ and $S$.  Lemma 4.2
then implies that the integral
$$
X\mapsto \int_{R_X} e^{(2\rho_+' - 2\rho_P)(Y)} dY, 
\quad X\in R_{i,j}(\delta', T-(T_2)_P^Q,S),
$$
is a fixed t-finite function of $X$, $T-(T_2)_P^Q$, and $S$, for all
well-situated $(T,S)$.

We can therefore write (4.22) as 
$$
\eqalignno{
&\int_{\exp \oR(T)} e^{-2\rho_P(H)p(\overline a))} v(H_P(\oa),T,S)
\int_{\omega_P} \sum_{\g\in\fo\cap W_0'(\bQ)} &(4.23)\cr
&\quad \sum_{\g_+\in U_0'(\bQ)} \int_{W_+(\bA) \oplus V_+(\bA)} 
f\bigl( \pi(n\oa)^\1 (\g + v) \bigr) \psi(\g_+ \cdot v)\,dv\,dn\,d\oa,\cr}
$$
for a fixed t-finite function $v$ on $\fa_P\times\fa\times\fa$.
This is still the integral of a continuous function on a compact set and so
converges absolutely.  At this point, we are almost done.  

Write $\exp\oR$ for the set $\{\oa \mid \overline{H_P(a)} \in \oR_{i,j}\}$,
where $\oR_{i,j}$ is as at the end of section 3.  The sum over all
geometric equivalence classes of the absolute value of
$$
\eqalignno{
&\int_{\exp(\oR_{i,j}\setminus\oR(T))} e^{-2\rho_P(H_P(\oa))}
v(H_P(\oa),T,S) 
\int_{\omega_P} \sum_{\g\in\fo\cap W_0'(\bQ)} &(4.24)\cr
&\quad \sum_{\g_+\in U_0'(\bQ)} \int_{W_+(\bA) \oplus V_+(\bA)} 
f\bigl( \pi(n\oa)^\1 (\g + v) \bigr) \psi(\g_+ \cdot v)\,dv\,dn\,d\oa\cr}
$$
converges absolutely and can be shown to be an error term, since every
point $X$ in $\oR\setminus\oR(T)$ satisfies (4.19), so that the argument
following (4.19) again applies.

However, the sum of (4.23) and (4.24) is the (absolutely convergent)
integral
$$
\eqalignno{
&\int_{\exp\oR_{i,j}} e^{-2\rho_P(H_P(a))} v(H_P(\oa),T,S) 
\int_{\omega_P} \sum_{\g\in\fo\cap W_0'(\bQ)} &(4.25)\cr
&\quad \sum_{\g_+\in U_0'(\bQ)} \int_{W_+(\bA) \oplus V_+(\bA)} 
f\bigl( \pi(n\oa)^\1 (\g + v) \bigr) \psi(\g_+ \cdot v)\,dv\,dn\,d\oa,\cr}
$$
which is a t-finite function in $T$ and $S$, since its dependence on them
arises only through the function $v(\cdot,T,S).$  The seminorm $\|\cdot\|$
needed in the statement of the theorem can be chosen to be the sum of all
the seminorms that appeared when bounding each error term.

The sum of (4.25) over all geometric equivalence classes converges, and is
again t-finite, as the function $v$ does not depend on $\fo$.  This
completes the proof of the Theorem. 
\qed\enddemo

\remark{Remark}  The proof of theorem implies (just as in [8]) that the
integral 
$$
\int_{\Gqa} \gv f(\pi(g^\1)\g)
$$
converges if and only if the linear functional
$$
\sum_{\lambda\in \Pi} \max (m_\lambda \lambda, 0) 
 - \sum_{\alpha\in\Sigma} \max (m_\alpha\alpha, 0) 
$$
is negative on $\fa\setminus \{\,0\,\}$, where $m_\lambda$ and $m_\alpha$
denote 
the multiplicity of the weight in the representations $\pi$ and $\Ad$,
respectively. The sufficiency of this condition for convergence of the
integral is due to Weil [13].  Its necessity was apparently also
known, and is due to Igusa%%% [?]
.
\endremark

Fix one $\pi$-dependent cone $\C.$  The above theorem shows that on this
cone, the functions $J_\fo^T(f,\pi)$ approximate t-finite functions.
As in [8], the proof of the theorem allows us to explicitly produce the
non-constant terms of each of the t-finite functions $P_{\fo,\C}$, so all 
we need to completely determine the functions $P_{\fo,\C}$ is the constant
term with respect to any point in $\fa.$  Let $T_0$ be the unique point in
$\fa$ such that 
$$
H(w_s^\1) + s^\1 T_0 = T_0
$$
for every element $s$ of the Weyl group of $(G,A)$, where $w_s$ is any
representative of $s$ in $G(\bQ)$; the existence of $T_0$ is the statement of
Lemma 1.1 of [2].  Write $P_{\fo,\C}(T)$ as a finite linear combination
of functions $e^{\lambda(T - T_0)} (T-T_0)^n,\ \lambda \in \fa^*,\ n$ a
nonnegative integer, and set $J_{\fo,\C}$ to be the constant term, that is the
term where both $\lambda$ and $n$ equal zero. Then the basic
form of the truncated Poisson summation formula for the representation
$\pi$ of $G$ on $V$ and the function $f$ on $V(\bA)$ is the following
theorem, proven exactly as in [8].

\proclaim{Theorem 4.3}
$$
\sum_{\fo \in \fo} J_{\fo,\C} (f,\pi ) \ = \ \sum_{\hfo \in \hfO} 
J_{\hfo,\C} (\hf , \vpi ). 
$$
\endproclaim

\remark{Remark}  Notice that because the weights of $\vpi$ are the negatives
of the weights of $\pi$, that the cones $\C$ determined by $\pi$ and $\vpi$
are the same. 
\endremark

The definition of $J_{\fo,\C} (f,\pi )$ depended on a number of choices.
The methods of [8] (based on those in [2]) show that the distributions
$J_{\fo,\C}$ are independent of $\omega$ and $T_1$, and that if $s$ is an
element of the Weyl group of $(G,A)$ and $J'$ denotes the
constant term of the truncated integral with respect to the non-standard
minimal parabolic subgroup $w_s^\1 P_0 w_s$, then
$$J_{\fo,\C} (f,\pi ) \ = \ J_{\fo,s^\1\C}' (f,\pi ).$$
If the representation $\pi$ is the Adjoint representation, then this
formula does depend on our choice of $K$, but for other
representations it need not.

\beginsection 5. Bibliography.

\ref \no 1 \by  J. Arthur \paper A trace formula for reductive groups I:
terms associated to classes in $G(\bQ)$ \jour Duke Math J. \vol 45 \yr1978
\pages 911--952 \endref

\ref \no 2 \by J. Arthur \paper The trace formula in invariant form \jour
Annals of Math. \vol 114 \yr 1981 \pages 1--74 \endref

\ref \no 3 \by J. Arthur \paper A measure on the unipotent variety \jour
Canad. J. Math \vol 37 \yr 1985 \pages 1237--1274\endref

\ref \no 4 \by M. Brion and M. Vergne \paper Residue formulae, vector
partition functions, and lattice points in rational polytopes \jour
J. Amer. Math. Soc. \vol 10 \yr 1997 \pages 797--833\endref

\ref \no 5 \by  A. Br\"ondsted \book An introduction to convex polytopes
\bookinfo Graduate Texts in Mathematics \vol 90 
\publ Springer-Verlag, New York-Berlin \yr 1983\endref
  
\ref \no 6 \by G. Kempf \paper Instability in invariant theory \jour Annals
of Math. \vol 108 \yr1978 \pages 299--316 \endref

\ref \no 7 \by S. Kudla and S. Rallis \paper A regularized Siegel-Weil
formula: The first term identity \jour Annals of Math. \vol 140 \yr 1994
\pages 1--80 \endref

\ref \no 8 \by J. Levy \paper A truncated Poisson formula for groups of
rank at most two \jour Amer. J. Math. \vol 117 \yr 1995 \pages
1371--1408 \endref 

\ref \no 9 \by J. Levy \paper Rationality of orbit closures\endref

\ref \no 10 \by D. Luna \paper Sur certaines op\'erations diffe\'rentiables
des
groupes de Lie \jour Amer. J. Math. \vol 97 \yr 1975 \pages 172--181
\endref

\ref \no 11 \by C. Rader and S. Rallis \paper Spherical characters on
$p$-adic symmetric spaces \jour Amer. J. Math. \vol 118 \yr 1996 \pages 
91--178 \endref 

\ref \no 12 \by R. W. Richardson \paper Conjugacy classes of $n$-tuples in
Lie algebras and algebraic groups \jour Duke Math. J. \vol 57 \yr 1988
\pages 1--35 \endref

\ref \no 13 \by A. Weil \paper Sur la Formule de Siegel dans la
Th\'eorie des Groupes Classiques \jour Acta Math. \vol 113 \yr 1965 \pages
1--87 \endref 

\ref \no 14 \by A. Yukie \book Shintani zeta functions \bookinfo London
Mathematical Society Lecture Note Series \vol 183 \publ Cambridge
University Press \yr 1993\endref

\enddocument